\PassOptionsToPackage{square,numbers}{natbib} 
\documentclass[10pt]{article}

\usepackage[utf8]{inputenc}
\usepackage[T1]{fontenc}
\usepackage[english]{babel}
\usepackage{graphicx}
\usepackage{amsmath,amssymb,amsbsy,amsthm,color,natbib}
\usepackage{multicol}
\usepackage{color} 
\usepackage{natbib}
\usepackage[linktoc=all,colorlinks=true,citecolor=blue,linkcolor=blue]{hyperref}
\usepackage[paperwidth=18cm,paperheight=29.7cm,hmargin=2cm,vmargin=3cm]{geometry}
\usepackage{authblk}
\usepackage{comment}
\usepackage{appendix}
\usepackage{bm}
\usepackage{makecell}

\graphicspath{{figuras/}}

\title{\textbf{Learning the Chaotic and Regular Nature \\ of Trajectories in Hamiltonian Systems with Lagrangian descriptors}}

\author[1]{Javier Jiménez López\thanks{javier.jimenezl@edu.uah.es}}
\author[1]{V. J. García-Garrido\thanks{ vjose.garcia@uah.es}}
\affil[1]{Departamento de F\'isica y Matem\'aticas, Facultad de Ciencias, Universidad de Alcal\'a, 28805 Alcal\'a de Henares, Madrid, Spain.}

\begin{document}

\maketitle


\begin{abstract}

In this paper, we explore the application of Machine Learning techniques, specifically Support Vector Machines (SVM), to unveil the chaotic and regular nature of trajectories in Hamiltonian systems using Lagrangian descriptors. Traditional chaos indicators, while effective, are computationally expensive and require an exhaustive study of the parameter space to establish the classification thresholds. By using SVMs trained on a dataset obtained from the analysis of the dynamics of the double pendulum Hamiltonian system, we aim at reducing the complexity of this process. Our trained SVM models demonstrate high accuracy when it comes to classifying trajectories in diverse Hamiltonian systems, such as for example in the four-well Hamiltonian, the Hénon-Heiles system and the Chirikov Standard Map. The results indicate that SVMs, when combined with Lagrangian descriptors, offer a robust and efficient method for chaos classification across different dynamical systems. Our approach not only simplifies the classification process but also is highlighting the potential of Machine Learning algorithms in the study of nonlinear dynamics and chaos. 

\end{abstract}

\noindent \textbf{keywords: Hamiltonian systems, chaos indicators, supported vector machines, Lagrangian descriptors}.


\section{Introduction} \label{Introduction}

The study of dynamical systems is a key area in contemporary physics due to its applications in other fields such as chemistry \cite{katsanikas2020detection}, life sciences \cite{coffey1998self, dokoumetzidis2001nonlinear}, or engineering \cite{parker1987chaos} as well as its contribution to the understanding of nature and its behavior. Given the fact that the vast majority of systems in nature are chaotic, it is of great interest to understand the origin of this behavior and how can we detect and measure it.

Until the present day, a wide range of tools known as chaos indicators have been developed \cite{skokos2016chaos}, allowing the study of the chaotic or regular nature of an initial condition for a given a dynamical system. These indicators are usually computationally expensive to calculate since it is necessary to conduct a study of the temporal evolution of the trajectory. Even for indicators derived from Lagrangian descriptors, which are known to be very efficient, evolving the trajectories to classify them as regular or chaotic is required. These indicators, although very efficient, have the problem of needing a threshold to classify the initial condition that is not simple to calculate as we need to have several initial conditions for the same energy and repeat this for a big set of energies in order to provide the algorithm that calculates the threshold with enough data to ensure a correct result, as we have already shown in \cite{jimenez2024pendulum}. This algorithm searches for the minimum value in the histogram generated when plotting the logarithm in base $10$ values of the corresponding chaos indicator for all the initial conditions simulated for an specific energy. For each energy value, the algorithm improves the precision of the threshold until it becomes close to constant, being this the reason why it need's a wide range of energies to provide the correct results for the threshold.

In order to address this problem, we propose the use of Machine Learning techniques, as this approach would allow us to develop a trained model that acts as a classifier without the need to manually determine a threshold to distinguish between the regular and chaotic behavior of trajectories. These techniques, specially neural networks \cite{chen2019symplectic, zhu2020deep, david2023symplectic, chen2022learning, zhang2021learning}, have started to become popular in the field of dynamical systems for different purposes. For example, we can highlight their application for solving systems of ordinary differential equations \cite{bakthavatchalam2022primer, chen2018neural}, and for studying Hamiltonian systems \cite{david2023symplectic, mattheakis2022hamiltonian}, as they offer certain advantages over traditional numerical methods. Other problems that convolutional neural networks have helped to address is the detection and classification of chaos \cite{Seok2020,celletti2022,Barrio2023,mayora2024} by means of using the Lyapunov exponents or other diagnostics such as the frequency map.  

In this work we have chosen to use Support Vector Machines (SVM), which have already been successfully applied in other nonlinear dynamics problems \cite{chan2001modelling, naik2021support,krajnak2022}, due to their high performance and conceptual simplicity. To train the model, we have used a varying number of classified initial conditions of the  double pendulum \cite{jimenez2024pendulum} and used it as a classifier for a system whose potential energy function has several critical points so that it can show diverse behaviors \cite{garcia2020exploring} depending on three know parameters of the model. In addition to the test on the above mentioned system, also known as the four-well Hamiltonian, we have also performed tests on the H\'enon-Heiles system \cite{henon1964applicability} and in the $2$D Standard Map \cite{Chirikov1979}, showing that SVMs provide great results even when compared with much more complex Machine Learning and Deep Learning techniques \cite{Barrio2023, Seok2020}. With this approach, we will be able to analyze if chaos can be characterized in the same way for  these systems or if, on the contrary, there is not any relation at all between the dynamical behavior they exhibit.

This paper is organized as follows. Section \ref{Methodology} is devoted to describing the methodology we have followed in this work. First we introduce the different chaos indicators used for this study, starting with the Smaller Alignment Index (SALI), and continuing with those based on the method of Lagrangian descriptors. In the second part of this section we show what a Support Vector Machine (SVM) is and how to implement it using neural networks in order to run the training procedure of the different models in parallel. In Sec. \ref{Results} we present and analyze the main results obtained regarding the chaotic and regular classification of trajectories for all the Hamiltonian systems we have considered. We finish with Sec. \ref{Conclusions}, where we provide a summary and outlook of the results obtained in this paper, and we discuss some ideas that we will explore in the future. The reader can find all the relevant information about the different Hamiltonian systems used in this work in Appendix \ref{Appendix}, and details about the simulations carried out for each system. 


\section{Methodology} \label{Methodology}

This section describes the methodology that we have followed to develop, train and test the SVM model with the capability of distinguishing between the regular and chaotic behavior of the trajectories of a given dynamical system. Since the SVM model and all the analysis carried out throughout this work uses data obtained from different chaos indicators, in particular the Smaller Alignment Index (SALI) and those derived from Lagrangian descriptors (LDs), we begin by explaining the basics of these diagnostic tools.  

\subsection{Chaos Indicators}

Identifying whether a deterministic dynamical system exhibits chaotic or regular behavior, and pinpointing areas within the phase space where instabilities are likely to emerge, is critically important across various disciplines, including astronomy, particle physics, and climate science. Chaos limits our capability to predict the system's future state accurately at various scales. In many real-world scenarios, understanding the impact of chaos on the system's overall dynamics is essential. Therefore, developing precise and efficient numerical tools to distinguish order from chaos, both on local and global scales, is crucial, especially for multidimensional systems with complex phase spaces. Today, numerous methods exist to address this challenge.

\subsubsection{The Smaller Alignment Index}

The Smaller Alignment Index (SALI) is a chaos indicator for Hamiltonian systems \cite{arnold2013mathematical} and symplectic maps \cite{Meiss1992} developed originally by Ch. Skokos \cite{Skokos2001}. This diagnostic tool is a reliable indicator because it provides a clear and distinct measure for differentiating between regular and chaotic motion of trajectories for any dynamical system \cite{skokos2016chaos}. For this reason, in this work, we regard the results it generates for the classification of trajectories as the ground truth to compare with the predictions made by our SVM models. The main idea behind SALI is very simple and is based on tracking the parallel or antiparallel alignment of two deviation vectors constructed for the reference trajectory under study. In the case that the trajectory is chaotic, the deviation vectors tend to align with the closest unstable manifold, while for a regular orbits, the deviation vectors eventually become tangent to the corresponding torus.

\begin{figure}[!h]
    \centering
    \includegraphics[scale = 0.45]{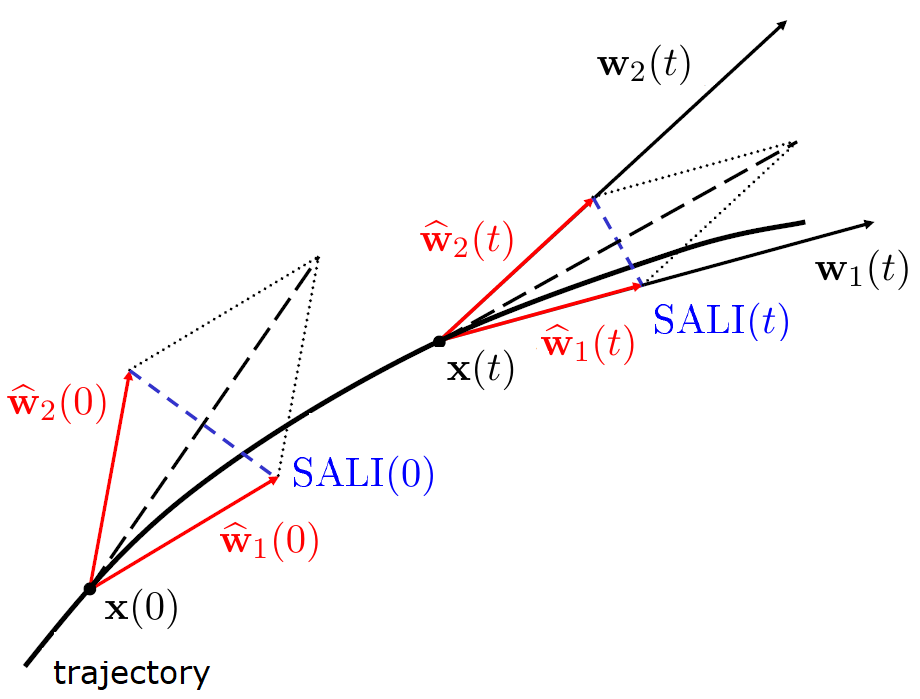}
    \caption{Evolution of a trajectory starting at the initial condition $\mathbf{x}(0)$, and also of two deviation vectors $\mathbf{w}_1(t)$ and $\mathbf{w}_2(t)$ that characterize the behavior of two neighboring trajectories. The SALI indicator at each time $t$ corresponds to the minimum length of the two diagonals of the parallelogram generated by the normalized deviation vectors. Figure adapted from \cite{skokos2016chaos}.}
    \label{SALI_diagram}
\end{figure}

The process for computing SALI to classify the chaotic or regular nature of an orbit is straightforward and consists in following simultaneously the time evolution of the trajectory starting at an initial condition $\mathbf{x}(0)$ (using the equations of motion that define the dynamical system), and of two initial deviation vectors $\mathbf{w}_1(0)$ and $\mathbf{w}_2(0)$ (for this task one needs to solve the variational equations) that we can consider initially to be orthonormal to each other. Since we are only interested in the direction of the deviation vectors along their evolution, we can normalize them at each time step, setting:
\begin{equation}
    \mathbf{\widehat{w}}_i(t) = \dfrac{\mathbf{w}_i(t)}{\lVert \mathbf{w}_i(t) \rVert} \;\;,\;\; i=1,2 \;,
\end{equation}
to control their exponential growth. If we introduce the parallel ($d_+$) and antiparallel ($d_{-}$) indices as follows:
\begin{equation}
    d_+(t) = \lVert \mathbf{\widehat{w}}_1(t) - \mathbf{\widehat{w}}_2(t) \rVert \quad,\quad d_-(t) = \lVert \mathbf{\widehat{w}}_1(t) + \mathbf{\widehat{w}_2}(t) \rVert \;,
\end{equation}
and define the SALI indicator as:
\begin{equation}
    \text{SALI}(t) = \min\left\lbrace d_+(t),d_-(t) \right\rbrace \;.
\end{equation}
A schematic representation of the procedure to compute SALI is given in Fig. \ref{SALI_diagram}.

It is a well known theoretical result \cite{skokos2016chaos} that, for two-dimensional symplectic maps, the SALI indicator displays the following asymptotic behavior:
\begin{equation}
    \text{SALI}(t) \; \propto \; \begin{cases}
        1/t^2 &,\;\; \text{for regular orbits} \\
        e^{-2\lambda_1 t} &,\;\; \text{for chaotic orbits} 
    \end{cases}
    \label{asymp_SALI_map}
\end{equation}
where $\lambda_1$ is the largest Lyapunov exponent of the orbit, while for continuous Hamiltonian systems it satisfies:
\begin{equation}
    \text{SALI}(t) \; \propto \; \begin{cases}
        \text{constant} &,\;\; \text{for regular trajectories} \\
        e^{-\lambda_1 t} \;&,\;\; \text{for chaotic trajectories} 
    \end{cases}
    \label{asymp_SALI_contHam}
\end{equation}
This asymptotic trends show that for a chaotic trajectory, SALI tends to zero, and this is a consequence of the fact that in this case the deviation vectors align with the closest unstable manifold. On the other hand, we see that for two-dimensional maps, the value of SALI for regular orbits also tends to zero, and this occurs because both deviation vectors align with the invariant torus, that in this case is a curve in the phase space. Finally, for a continuous Hamiltonian system, the value of SALI tends to a constant for regular trajectories. Therefore, the long-term behavior of SALI allows us to easily identify a trajectory as chaotic or regular, since the time evolution of the indicator occurs at different timescales. For this reason it is natural to employ the logarithm of SALI values to carry out the classification of trajectories.   

In order to illustrate how SALI can be used to characterize the chaotic or regular nature of trajectories, we present next a simple analysis for the Chirikov Standard Map (see the dynamical system described by Eq. \eqref{std_map} in the Appendix) and also for the H\'en{o}n-Heiles Hamiltonian (see Eq. \eqref{Ham_HH} in the Appendix). First we depict in Fig. \ref{PSOS_StdMap_K1p5} a Poincar\'e map of the Standard Map for $K = 1.5$, where we have marked a regular (blue dot) and a chaotic (orange dot) initial condition. In Fig. \ref{SALI_StdMap} A) we display the time evolution of $\log_{10}(\text{SALI})$ for the regular initial condition, and highlight its long-term behavior as described in Eq. \eqref{asymp_SALI_map}. We do the same in Fig. \ref{SALI_StdMap} B) but for the chaotic initial condition. Additionally, in panel C) of Fig. \ref{SALI_StdMap} we display the histograms of $\log_{10}(\text{SALI})$ values, corresponding to three different ensembles of $10^4$ random initial conditions generated for the values of the model parameter $K=0.5$, $K=0.971635$ and $K=1.5$. The dashed vertical line indicates the threshold used to distinguish between chaotic and regular orbits, which we have set for this system to be $\log_{10}(\text{SALI}) = -13$. 

\begin{figure}[!h]
    \centering
    \includegraphics[scale = 0.29]{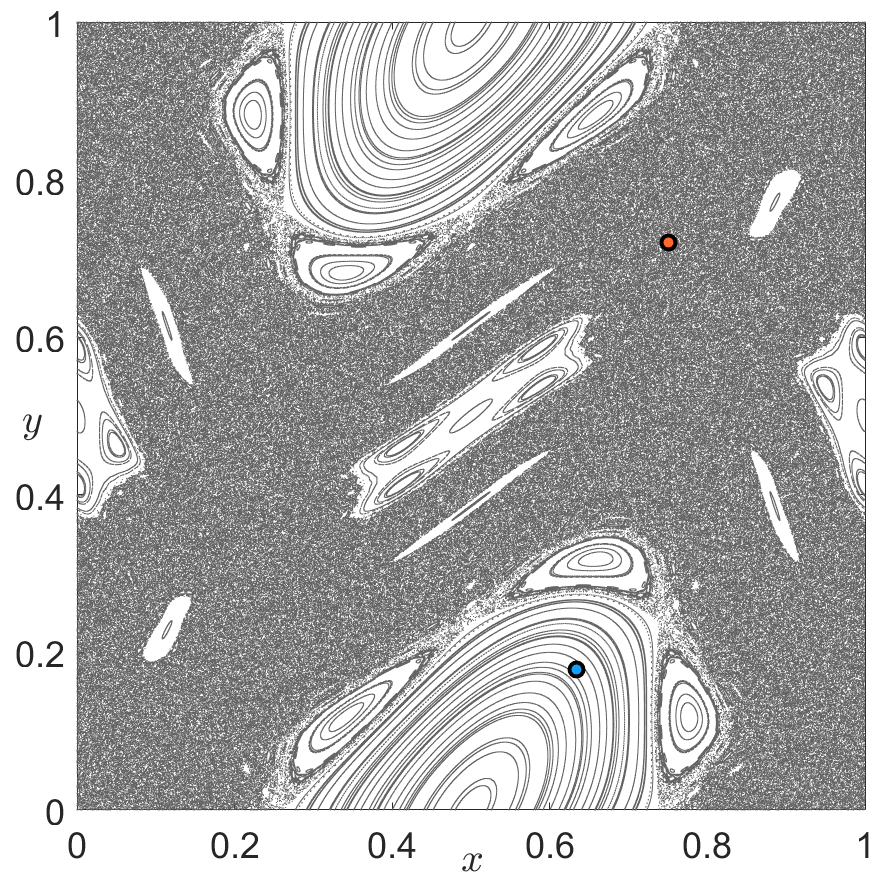}
    \caption{Poincar\'e map for the Standard Map in Eq. \eqref{std_map} with $K=1.5$. We have marked in the plot a regular (blue dot) and a chaotic initial condition (orange).}
    \label{PSOS_StdMap_K1p5}
\end{figure}

\begin{figure}[!h]
    \centering
    A) \includegraphics[scale = 0.25]{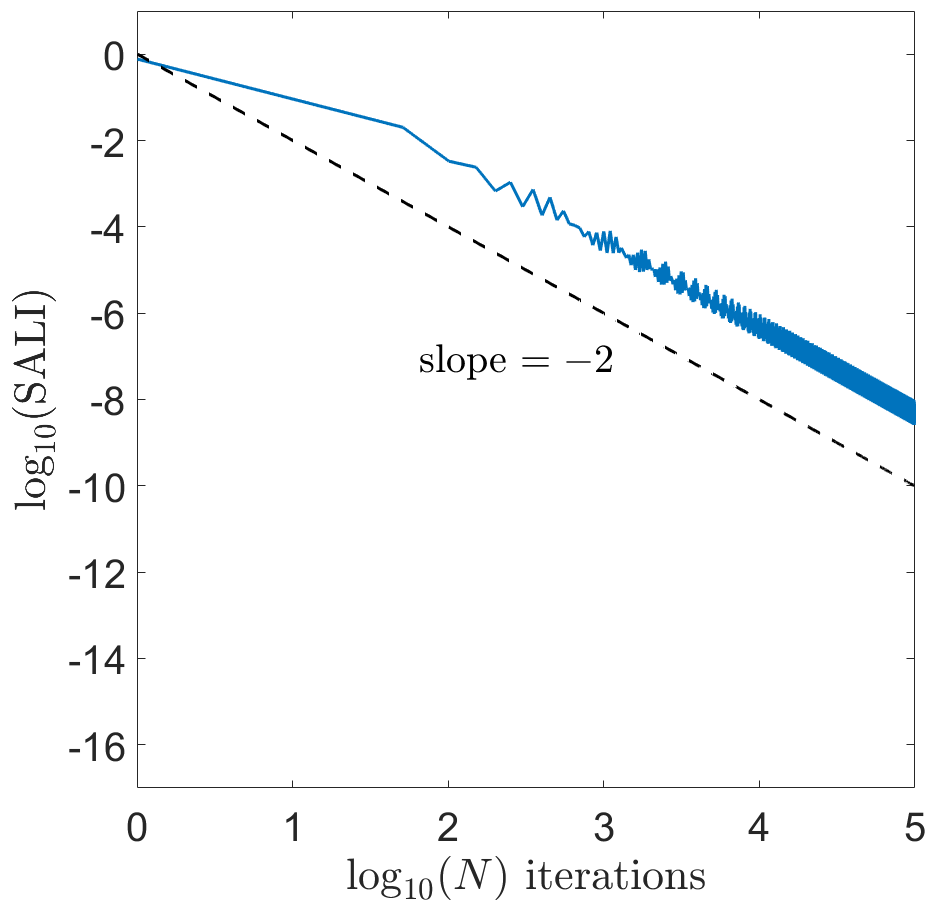}
    B) \includegraphics[scale = 0.25]{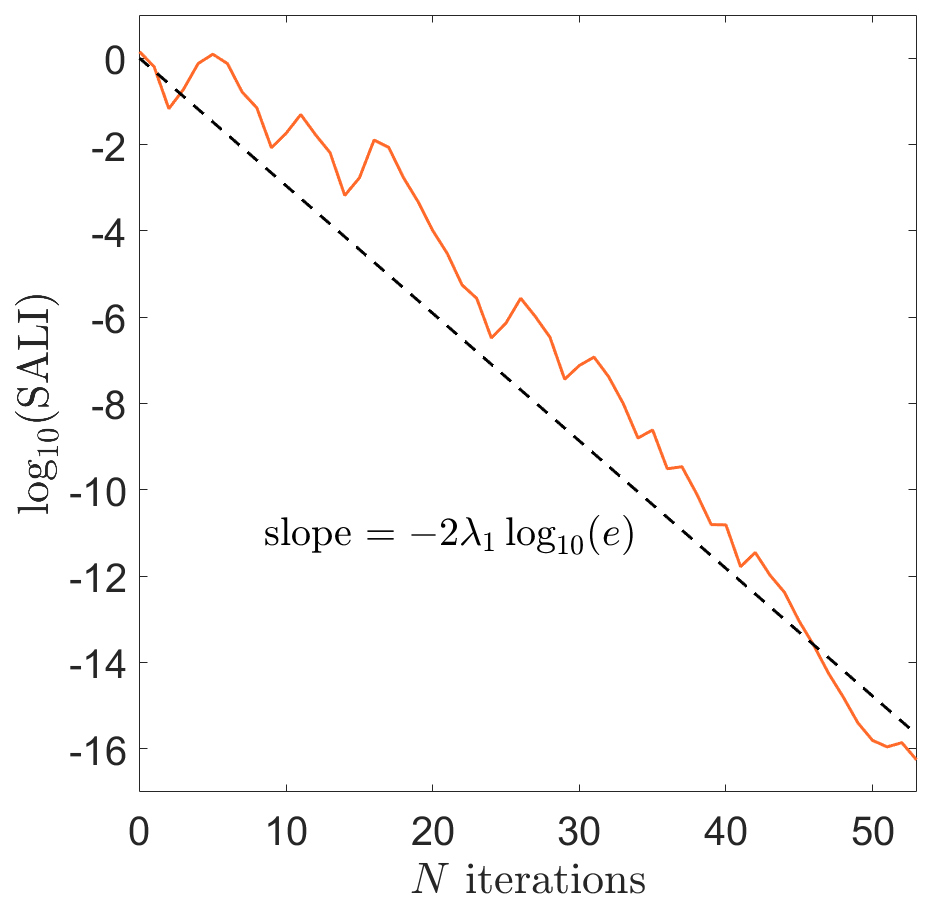} \\
    C) \includegraphics[scale = 0.28]{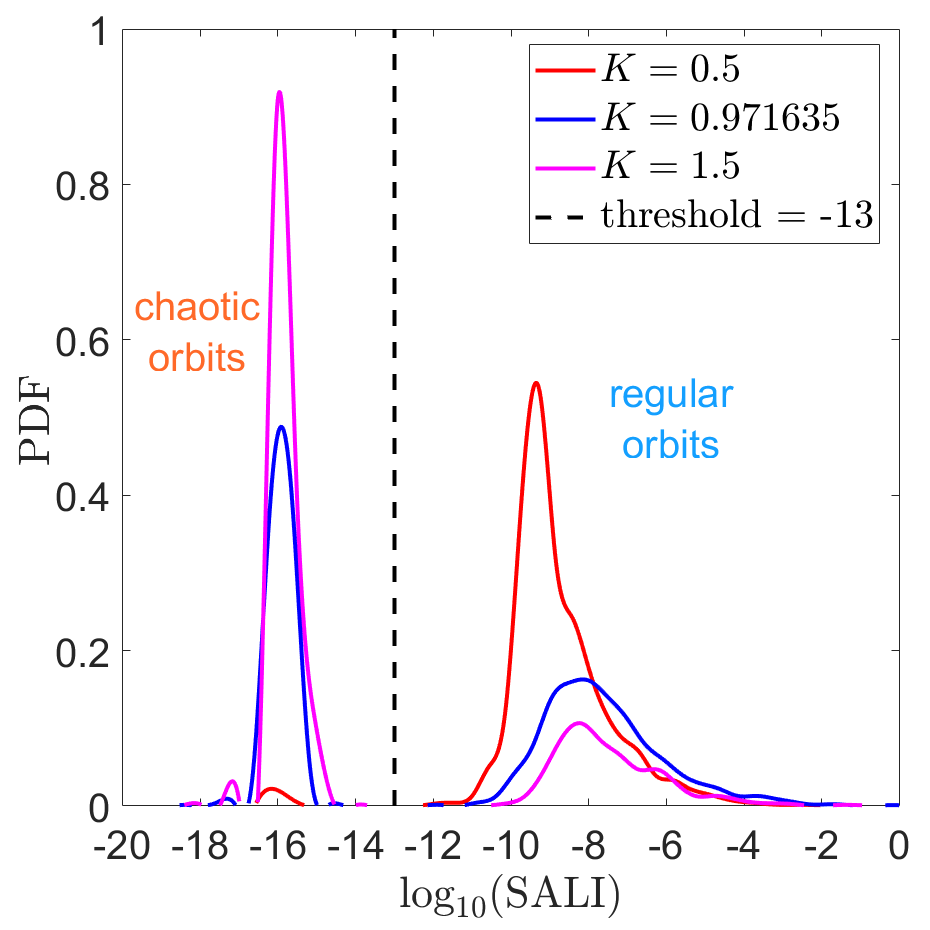}
\caption{A) Time evolution of $\log_{10}(\text{SALI})$ for the regular initial condition in Fig. \ref{PSOS_StdMap_K1p5}. The dashed line highlights the asymptotic behavior given in Eq. \eqref{asymp_SALI_map}; B) Time evolution of $\log_{10}(\text{SALI})$ for the chaotic initial condition in Fig. \ref{PSOS_StdMap_K1p5}. The dashed line shows the asymptotic behavior described in Eq. \eqref{asymp_SALI_map}; C) Histogram of $\log_{10}(\text{SALI})$ values, using a random ensemble of $10^4$ initial conditions for each of the cases simulated for the Standard Map. The dashed line indicates the threshold used to distinguish between chaotic and regular orbits.}
    \label{SALI_StdMap}
\end{figure}

\begin{figure}[!h]
    \centering
    \includegraphics[scale = 0.28]{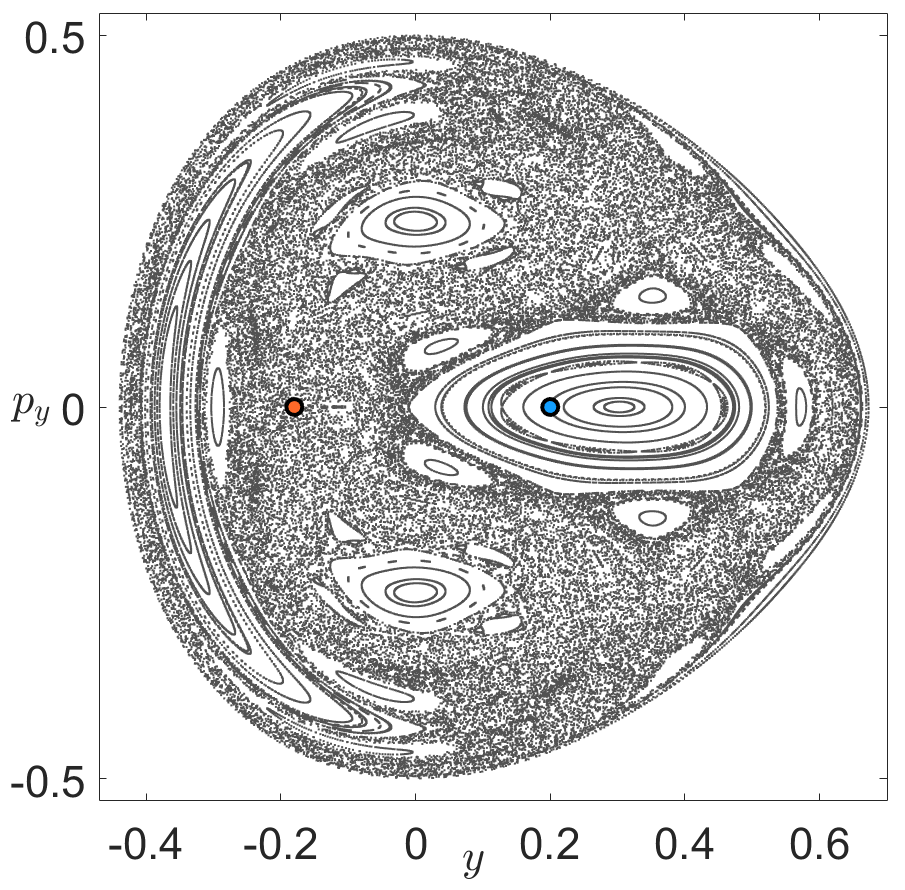}
    \caption{Poincaré map for the H\'{e}non-Heiles Hamiltonian in Eq. \eqref{Ham_HH} calculated on the surface of section $x = 0$, $p_x \geq 0$. using an energy $\mathcal{H} = 1/8$. We have marked in the plot a regular (blue dot) and a chaotic initial condition (orange).}
    \label{HH_PSOS}
\end{figure}

\begin{figure}[!h]
    \centering
    A) \includegraphics[scale = 0.25]{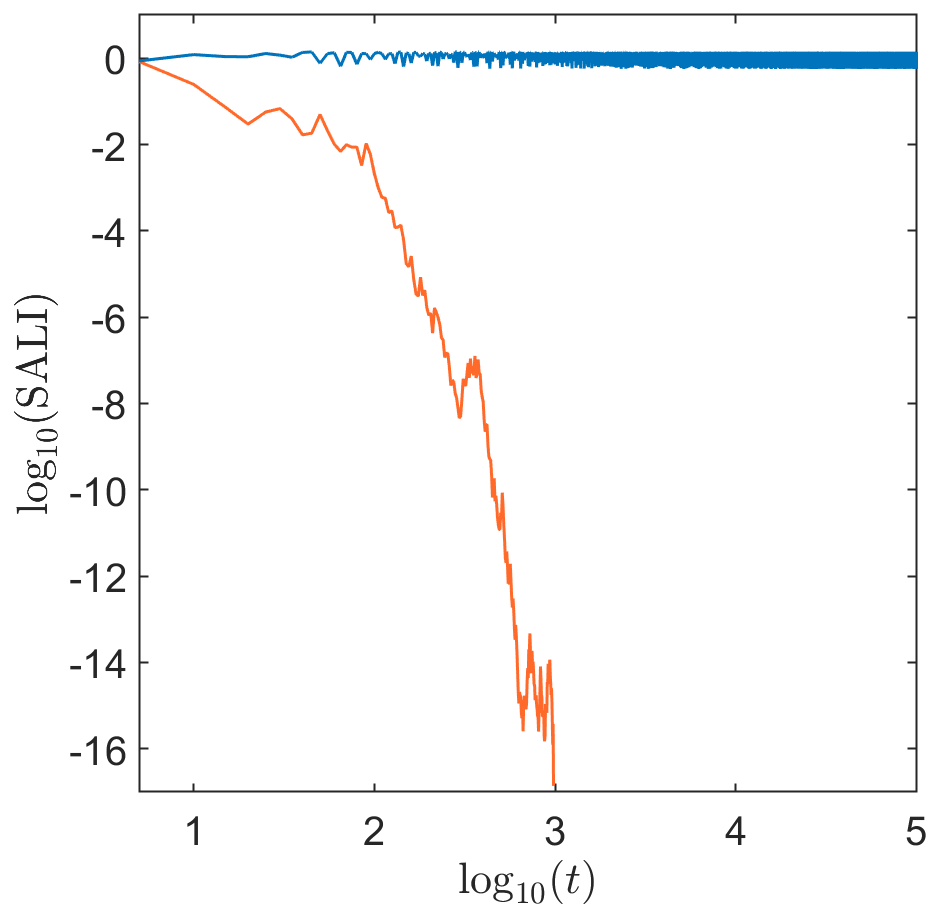}
    B) \includegraphics[scale = 0.25]{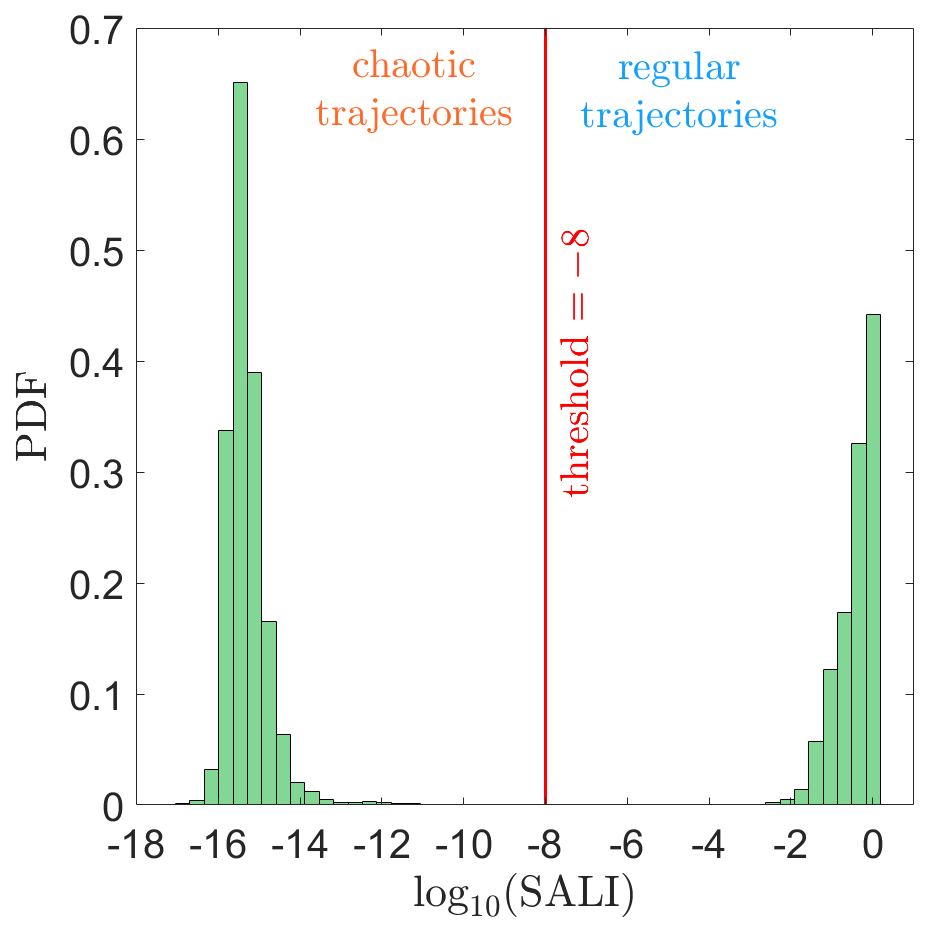}
    \caption{A) Time evolution of $\log_{10}(\text{SALI})$ for the regular (blue) and chaotic (orange) initial conditions depicted in Fig. \ref{HH_PSOS}; B) Histogram of $\log_{10}(\text{SALI})$ values, calculated for a random ensemble of $10^4$ initial conditions with energy $\mathcal{H} = 1/8$, where the threshold to differentiate chaotic and regular behavior is marked with a vertical dashed red line.}
    \label{SALI_HH_0p125}
\end{figure}

For the sake of completeness, we also demonstrate the behavior of SALI for regular and chaotic initial conditions in a continuous dynamical system, the classical H\'{e}non-Heiles Hamiltonian. In Fig. \ref{HH_PSOS} we display a Poincar\'e map on the surface of section $x = 0$, $p_x \geq 0$, calculated for this system using an energy $\mathcal{H}=1/8$, and we have superimposed a regular (blue dot) and chaotic (orange dot) initial condition. For these initial conditions we display in Fig. \ref{SALI_HH_0p125} A) the time evolution of SALI, where we can clearly observe the asymptotic behavior of SALI as described in Eq. \eqref{asymp_SALI_contHam}. Moreover, panel B) of Fig. \ref{SALI_HH_0p125} shows a histogram calculated with the values of $\log_{10}(\text{SALI})$ for a random ensemble of $10^4$ initial conditions selected in the same surface of section we have just described. Notice that in this case, we can take as a threshold to separate chaotic from regular trajectories the value $\log_{10}(\text{SALI}) = -8$, as it clearly divides the histogram into two disconnected parts, the one on the left corresponding to chaotic trajectories (those for which SALI exponentially converges to zero) and the one on the right for regular trajectories (those whose SALI value tends to a constant).   

\subsubsection{Lagrangian descriptors} \label{LD_ind}

\begin{figure}[!h]
    \centering
    \includegraphics[scale = 0.32]{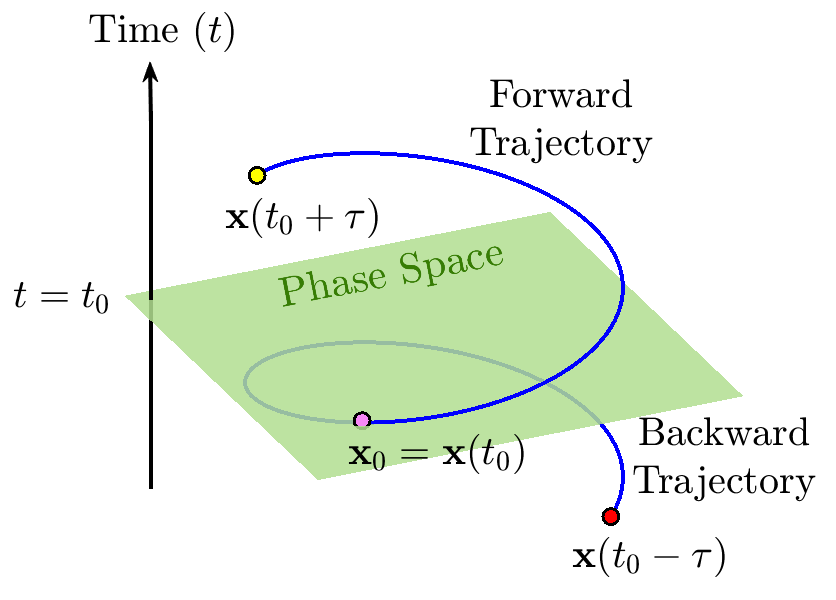}
    \caption{Initial condition $\mathbf{x}_0 = \mathbf{x}(t_0)$ at time $t = t_0$ evolving forward and backward in time for a given integration time $\tau$. A value for the Lagrangian descriptor function is assigned to the phase space point $\mathbf{x}_0$ by accumulating a non-negative function $\mathcal{F}(\mathbf{x},t)$ along its trajectory.}
    \label{LD_diagram}
\end{figure}

The method of Lagrangian descriptors (LDs) is a trajectory-based diagnostic technique that was originally developed in the field of Geophysics to analyze Lagrangian transport and mixing processes in the ocean and the atmosphere \cite{madrid2009,mendoza10}. A Lagrangian descriptor is a scalar function constructed in the following way. Given a continuous dynamical system:
\begin{equation}
    \dfrac{d\mathbf{x}}{dt} = \mathbf{f}(\mathbf{x},t) \;.
    \label{general_ds}
\end{equation}
we define a non-negative function $\mathcal{F}(\mathbf{x}(t;\mathbf{x}_0),t)$ that depends on the initial condition $\mathbf{x}_0$ at time $t=t_0$. To determine the LD scalar field, denoted by $\mathcal{L}$, we set an integration time $\tau > 0$ and calculate:
\begin{equation}
    \mathcal{L}(\mathbf{x}_0,t_0,\tau) = \mathcal{L}_f(\mathbf{x}_0,t_0,\tau) + \mathcal{L}_b(\mathbf{x}_0,t_0,\tau) \;,
    \label{ld_def}
\end{equation}
where the forward ($\mathcal{L}_f$) and backward ($\mathcal{L}_b$) components of the LD function are given by:
\begin{equation}
\mathcal{L}_f(\mathbf{x}_0,t_0,\tau) = \int_{t_0}^{t_0+\tau} \mathcal{F}(\mathbf{x}(t;\mathbf{x}_0),t) \, dt \quad,\quad
\mathcal{L}_b(\mathbf{x}_0,t_0,\tau) = \int_{t_0-\tau}^{t_0} \mathcal{F}(\mathbf{x}(t;\mathbf{x}_0),t) \, dt \;. 
\end{equation}
As Eq. \eqref{ld_def} shows, the calculation of LDs involves the accumulation of the values taken by the function $\mathcal{F}$ along the trajectory starting at $\mathbf{x}_0$, as it evolves forward and backward in time (see Fig. \ref{LD_diagram}). In the literature it has been shown rigorously that the scalar field generated by this method has the capability of identifying the invariant sets (equilibria, stable and unstable manifolds, tori, periodic orbits, etc,) that characterize the dynamical behavior of trajectories in the phase space of the system \cite{mancho2013,lopesino2017}. 

In this work we will use the definition of LDs that is based on the $p$-norm of Functional Analysis \cite{lopesino2017}. Thus, we select as the non-negative function the expression:
\begin{equation}
    \mathcal{F}(\mathbf{x},t) = \sum_{i=1}^n |f_i(\mathbf{x},t)|^p \;\;,\;\; 0 < p \leq 1 \;,
\end{equation}
where $f_i$ is the $i$-th component of the vector field that determines the dynamical system in Eq. \ref{general_ds}. In particular we will set $p = 1/2$, which is commonly used because it yields a scalar field where the phase space structures (stable and unstable manifolds of normally hyperbolic invariant manifolds, KAM tori, etc.) are nicely detected.

Similarly, if we are working with a $d$-dimensional discrete dynamical system (a map) of the form:
\begin{equation}
\mathbf{x}_{n+1} = \mathbf{f}(\mathbf{x}_{n}) \;,\quad n=0,1,\ldots	\;,
\label{discrete_DS}
\end{equation}
and we fix a number of iterations $N>0$ both forward (using $\mathbf{f}$) and backward in time (applying $\mathbf{f}^{-1}$), we can define the discrete version of LDs \cite{Lopesino2015} as follows:
\begin{equation}
\mathcal{L}\left(\mathbf{x}_0,N\right) = \sum_{n=-N}^{N-1} \sum_{j = 1}^{d} |x^j_{n+1}-x^j_{n}|^{p}
\label{DLD}
\end{equation}
where $0<p\leq 1$. We will use in this work the case $p=1/2$ to construct our discrete LD function. Note that Eq. \eqref{DLD} can be split into two different terms:
\begin{equation}
\mathcal{L}\left(\mathbf{x}_0,N\right) = \mathcal{L}_{f}\left(\mathbf{x}_0,N\right) + \mathcal{L}_{b}\left(\mathbf{x}_0,N\right) \;, 
\end{equation}
where $\mathcal{L}_{f}$ and $\mathcal{L}_{b}$ quantify, respectively, the contributions to the LDs of the forward and backward iterations of the orbit starting at the initial condition $\mathbf{x}_0$. This yields:
\begin{equation}
\mathcal{L}_f = \sum_{n=0}^{N-1} \sum_{j = 1}^{d} |x^j_{n+1}-x^j_{n}|^{p} \quad,\quad \mathcal{L}_{b} = \sum_{n=-N}^{-1} \sum_{j = 1}^{d} |x^j_{n+1}-x^j_{n}|^{p} \;.
\end{equation}

Chaos indicators based on LDs with the capability of classifying the regular or chaotic nature of trajectories were developed in \cite{Hille22,zimper23}. These diagnostics have been benchmarked and validated against SALI, where they have shown to have a success rate higher than $90\%$ when compared to SALI. Among the advantages that this methodology brings, we can highlight that implementation of LDs is straightforward, since the value of the LD function can be directly calculated by adding an extra differential equation (for continuous systems) or difference equation (for maps) to the equations of motion that define the dynamical system under study. This offers an edge with respect to computational time, since the chaos indicators based on LDs do not require the computation of the time evolution of the deviation vectors (which involves the numerical solution of the variational equations) that other classical chaos indicators such as SALI need. This largely simplifies the complexity of the simulations and reduces the required CPU time.

In order to construct the chaos indicators from LDs, we need to obtain the neighbors of the initial condition that we would like to analyze. These neighbors are given by:
\begin{equation}
\mathbf{y}_i^{\pm} = \mathbf{x}_0 \pm \sigma_i \, \mathbf{e}_i \;\;,\;\; i=1,\ldots,n \;,
\end{equation}
where $\mathbf{e}_i$ is the $i$-th canonical basis vector in $\mathbb{R}^n$, and $\sigma_i$ represents the distance between the central point $\mathbf{x}_0$ (the initial condition we would like to classify) and its neighbors on the grid. The value of $n$ represents the dimension of the space where the initial conditions are selected. In this paper, since we are working with two-dimensional phase space slices, then $n = 2$, and thus each initial condition on the grid has 4 neighbors. Moreover, for our simulations we have chosen a value of $\sigma_i = 10^{-4}$, which ensures that the LD-based chaos indicators introduced below in Eq. \eqref{ld_chaos_inds} maximize their accuracy for correctly identifying chaos and regularity when compared with SALI. An extensive study of how the spacing between the initial condition and its neighbors on the grid influence the accurate classification of chaotic and regular trajectories can be found in \cite{Hille22,zimper23}. Using these points, we can construct the following chaos indicators based on LDs:
\begin{equation}
\begin{split}
    D_L^n(\mathbf{x}_0) &= \dfrac{1}{2n \mathcal{L}_{f}(\mathbf{x}_0)} \sum_{i=1}^{n} \big|\mathcal{L}_{f}(\mathbf{x}_0)-\mathcal{L}_{f}\left(\mathbf{y}_i^{+}\right)\big|+\big|\mathcal{L}_{f}(\mathbf{x}_0)-\mathcal{L}_{f}\left(\mathbf{y}_i^{-}\right)\big| \;, \\  
    R_L^n(\mathbf{x}_0) &= \Bigg|1-\dfrac{1}{2n \mathcal{L}_{f}(\mathbf{x}_0)} \sum_{i=1}^{n} \mathcal{L}_{f}\left(\mathbf{y}_i^{+}\right)+\mathcal{L}_{f}\left(\mathbf{y}_i^{-}\right) \Bigg| \;, \\[.1cm]
    C_L^n(\mathbf{x}_0) &= \dfrac{1}{2n} \sum_{i=1}^{n} \dfrac{\big|\mathcal{L}_{f}\left(\mathbf{y}_i^{+}\right)-\mathcal{L}_{f}\left(\mathbf{y}_i^{-}\right)\big|}{\sigma_i} \;, \\[.1cm]
    S_L^n(\mathbf{x}_0) &= \dfrac{1}{n} \sum_{i=1}^{n} \dfrac{\big|\mathcal{L}_{f}\left(\mathbf{y}_i^{+}\right)-2\mathcal{L}_{f}(\mathbf{x}_0)+\mathcal{L}_{f}\left(\mathbf{y}_i^{-}\right)\big|}{\sigma_i^2} 
\end{split}
\label{ld_chaos_inds}
\end{equation}
where $\mathcal{L}_{f}(\cdot)$ is the forward LD value calculated for an integration time $\tau$. Note that the chaotic or regular nature of a trajectory is equivalently characterized if we integrate forward or backward in time, so only one of the components of the LD function is required to determine the chaos indicators. For this work we will mainly use the $S^n_{L}$ chaos indicator, as it is the one with the highest success rate (above $92\%$) when compared with the classification of trajectories provided by SALI \cite{zimper23}. 

\begin{figure}[!h]
    \centering
    A) \includegraphics[scale = 0.29]{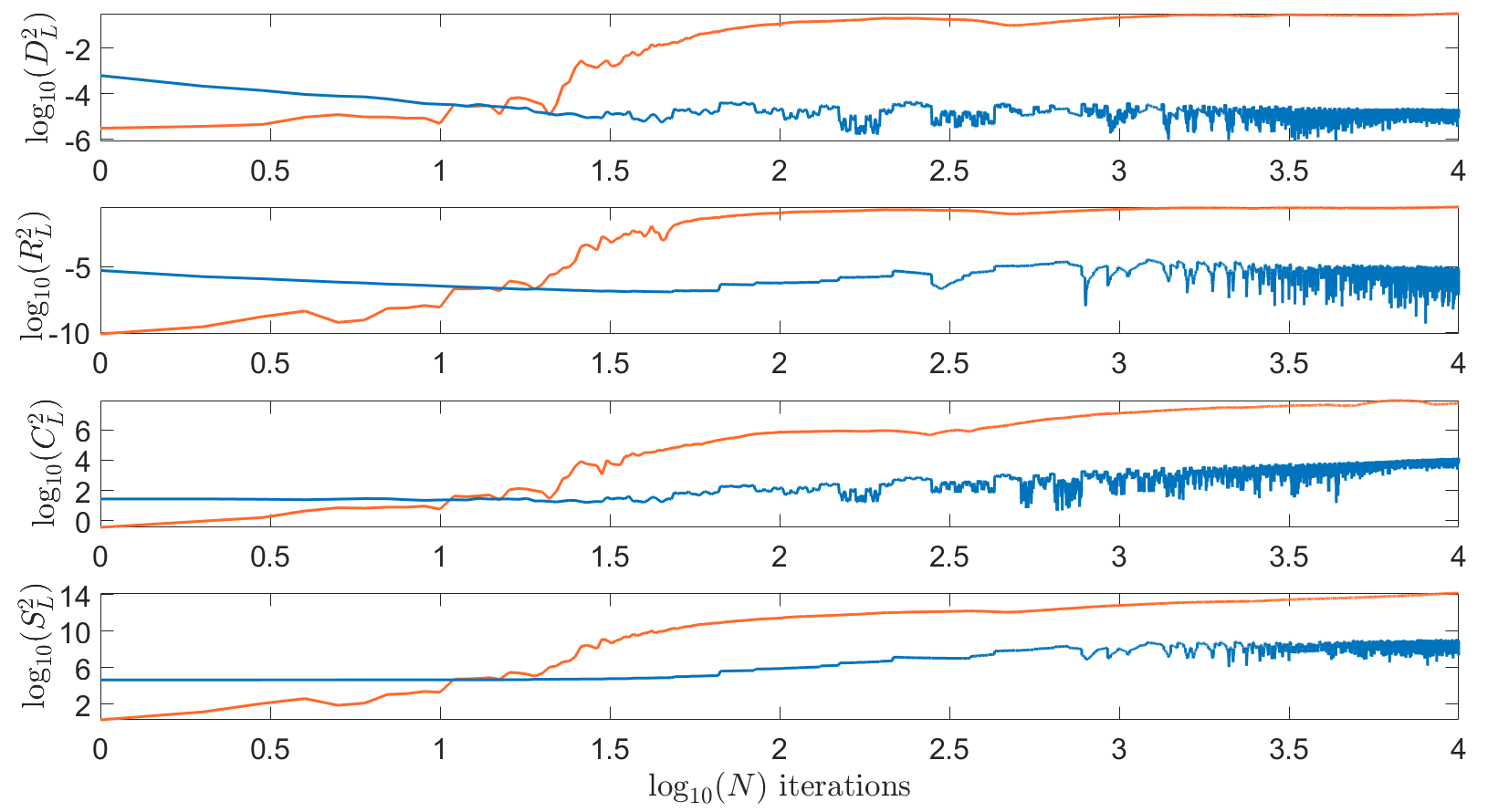} \\
    B) \includegraphics[scale = 0.28]{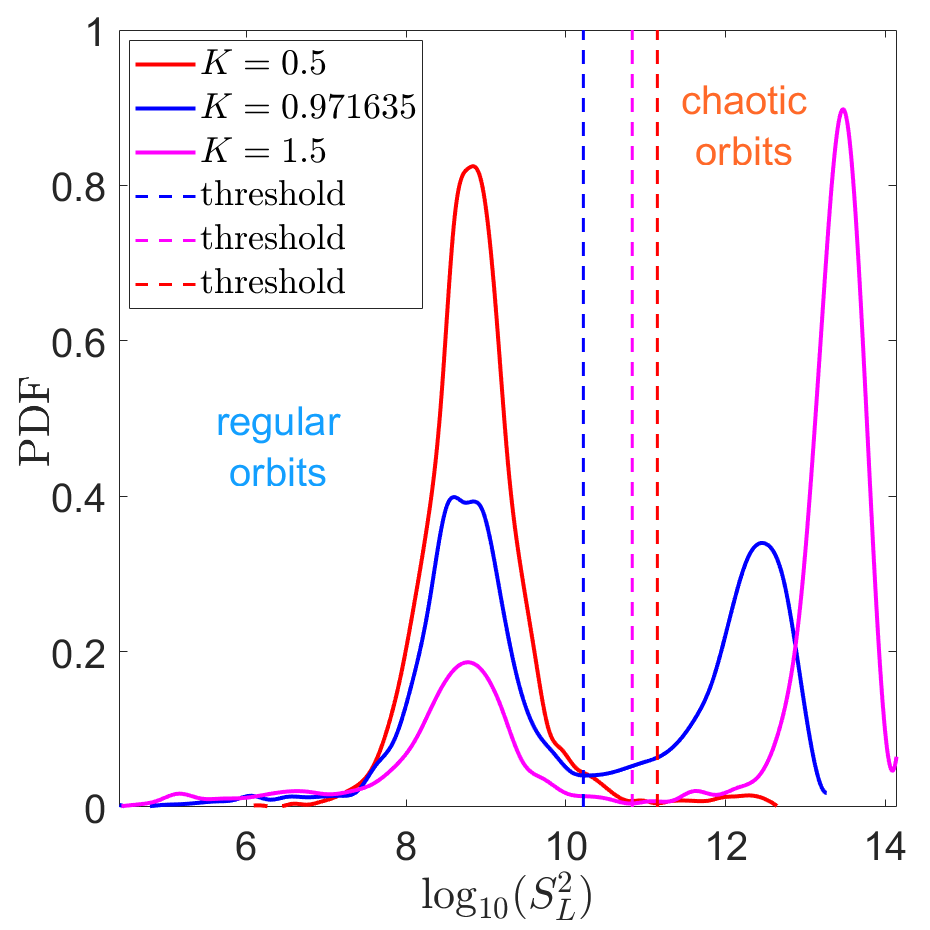}
    \caption{A) Time evolution of the logarithm of the LD-based chaos indicators defined in Eq. \eqref{ld_chaos_inds} calculated for the regular (blue) and chaotic (orange) initial conditions marked in Fig. \ref{PSOS_StdMap_K1p5} for the Standard Map with $K=1.5$; B) Histogram of $\log_{10}(S^2_L)$ values calculated for the Standard Map using three different random ensembles of $10^4$ initial conditions, sampled for the model parameter cases $K=0.5$, $K=0.971635$ and $K=1.5$. We have marked with vertical dashed lines the threshold values that separate chaotic from regular orbits in the system.}
    \label{Map}
\end{figure}

In order to illustrate how these chaos indicators can be used to distinguish between chaotic and regular behavior, we show in Fig. \ref{Map} A) their time evolution for a regular and a chaotic initial condition selected for the Standard Map in the case displayed in Fig. \ref{PSOS_StdMap_K1p5}. It is important to remark that the values taken by the indicators for the chaotic and regular initial conditions separate and become distinct with time. This important feature of the diagnostics is the one we will exploit in this work to develop an SVM model to accurately characterize the nature of the trajectories of a dynamical system. In Fig. \ref{Map} B) we depict the histogram of $\log_{10}(S_L^2)$ values, calculated for the Standard Map using three different random ensembles of $10^4$ initial conditions, sampled for the model parameter cases $K=0.5$, $K=0.971635$ and $K=1.5$. We can clearly see that the histogram has two peaks, which define the regular and chaotic orbits of the system. Hence, if we are interested in distinguishing order from chaos, we need to set a threshold value, that is taken in the histogram plot at the local minimum between both peaks. Note that this threshold depends not only on the dynamical system under study and on its parameters, but also on the chaos indicator chosen. By visual inspection of the histogram, it is not difficult to establish this threshold, however one needs a histogram generated by large ensembles of initial conditions to carry out this task. The main goal that we pursue in this paper is to build a SVM model that can identify chaotic and regular trajectories without establishing this threshold manually. 

\begin{figure}[!h]
    \centering
    A) \includegraphics[scale = 0.29]{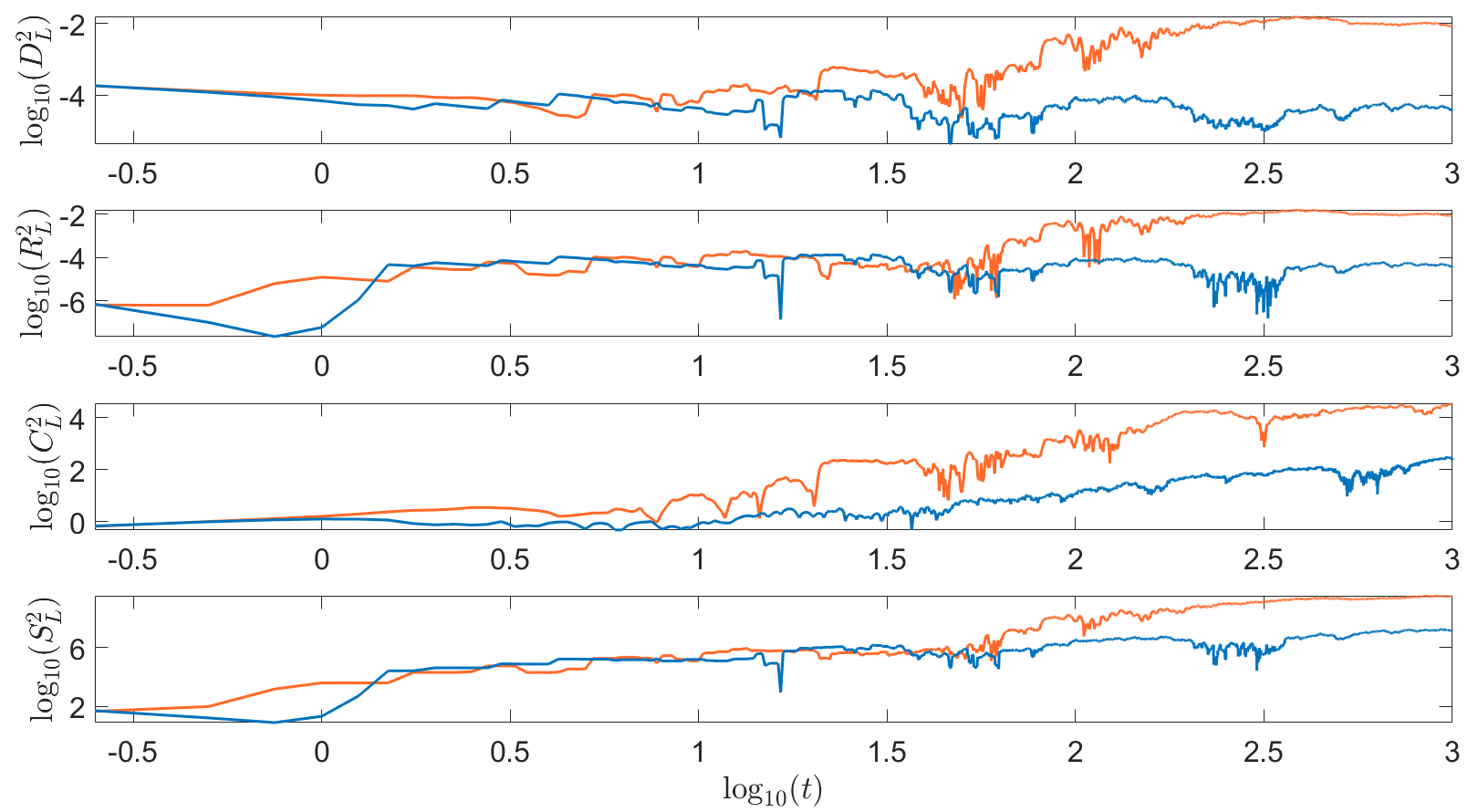} \\
    B) \includegraphics[scale = 0.28]{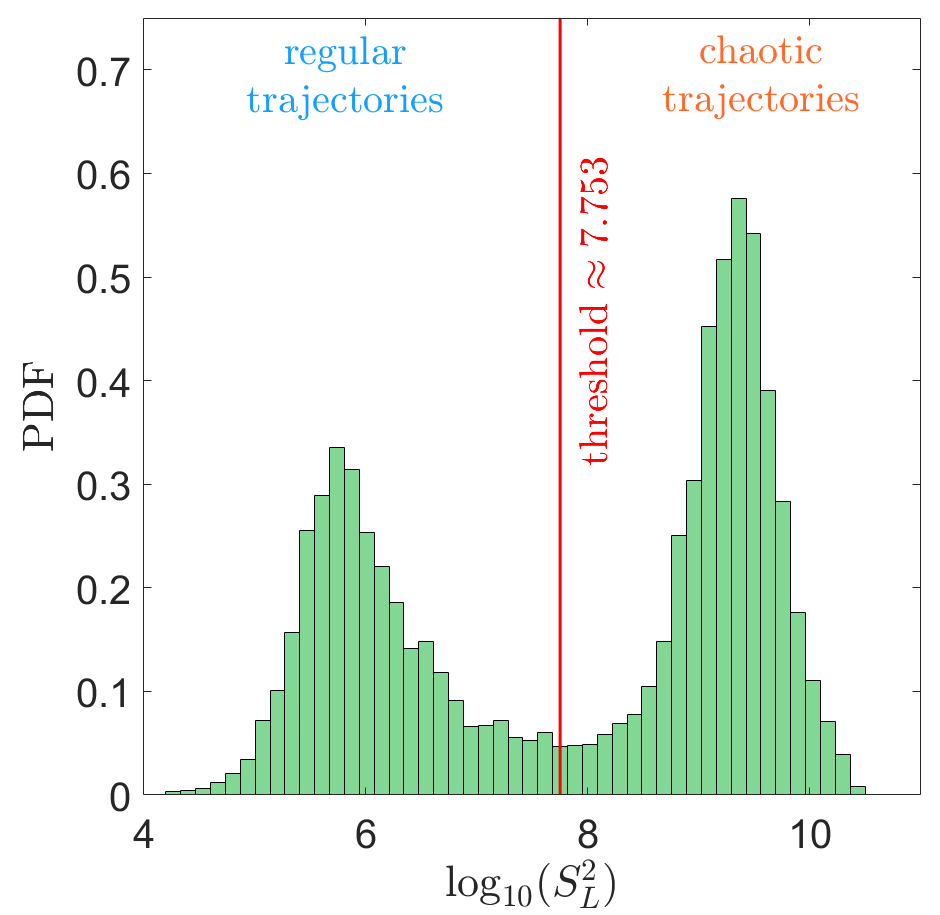}
    \caption{A) Time evolution of the logarithm of the LD-based chaos indicators defined in Eq. \eqref{ld_chaos_inds} calculated for the regular (blue) and chaotic (orange) initial conditions marked in Fig. \ref{HH_PSOS} for the H\'enon-Heiles Hamiltonian with energy $\mathcal{H}=1/8$; B) Histogram of $\log_{10}(S^2_L)$ values calculated for the H\'enon-Heiles Hamiltonian for a random ensemble of $10^4$ initial conditions, sampled for the energy $\mathcal{H}=1/8$ on the surface of section $x = 0$, $p_x \geq 0$. We have marked with a vertical red line the threshold value that separates chaotic from regular trajectories in the system.}
    \label{Henon}
\end{figure}

We finish this part of the methodology dedicated to explaining the basics of chaos indicator techniques by demonstrating in Fig. \ref{Henon} how the LD-based chaos indicators in Eq. \eqref{ld_chaos_inds} can also be easily applied to classify chaos and regularity in continuous systems, in particular we carry out this analysis for the H\'enon-Heiles Hamiltonian with energy $\mathcal{H} = 1/8$. 

\subsection{Machine Learning and Supported Vector Machines}

\begin{figure}[htbp]
	\centering
	\includegraphics[scale = 0.55]{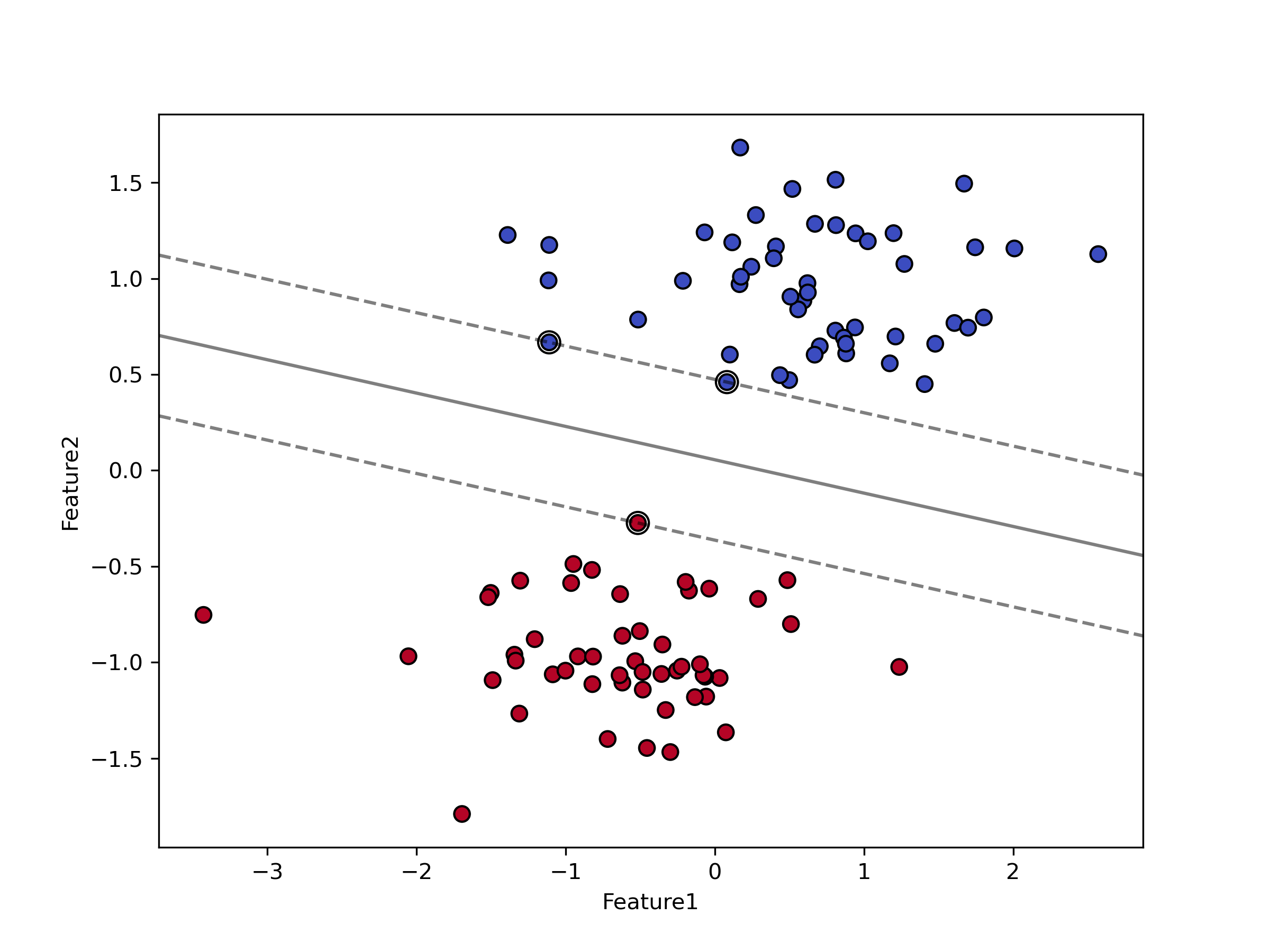}
	\caption{Example of how a SVM works with two different classes of data (blue and red) according to two features of the data (feature1 and feature2). In the graph, the decision plane and the support vectors are also showed, being the decision line the continuous one and the support vectors the dashed ones.}
	\label{SVM_Example}
\end{figure}

A Support Vector Machine (SVM) \cite{cortes1995support,mammone2009support,bishop2006pattern, rojo2018digital} is a widely-used supervised machine learning algorithm, particularly well-suited for binary classification tasks due to its simplicity and efficiency. The way this algorithm works is quite intuitive and its objective is to find a decision surface, usually a hyperplane even though more complex surfaces can also be used, that best separates the different types of data leaving the maximum margin between the classes. Consider a dataset $(x_{1},y_{1})$,$(x_{2},y_{2})$,$\ldots$,$(x_{n},y_{n})$ where $x_{i} \in \mathbb{R}^d$ represents the feature vectors and $y_{i} \in \{-1, 1\}$ are the binary class labels. What the algorithm will do is to find a function that assigns one of the two categories to new examples based on which side of the decision boundary the fall on. An example of it can be shown in Fig.\ref{SVM_Example}, where two types of data are classified according to this algorithm.

The mathematical foundation of SVMs is based on the kernel functions $\mathcal{K}(x,x^{\prime})$, which are a crucial concept in the theory of SVMs, allowing the them to efficiently handle both, linear and non-linear classification problems. These kernel functions compute the inner product of two vectors $x$ and $x^{\prime}$ in a high-dimensional space without explicitly mapping these vectors to that space. The specific mapping is done as follows:
\begin{equation} \label{Kernel}
    \mathcal{K}(x_{i}, x_{j}) = \langle \phi(x_{i}), \hspace{0.05cm} \phi(x_{j}) \rangle
\end{equation}
where $\phi : \mathbb{R}^d \longrightarrow \mathbb{R}^D$ with $D > d$. The most commonly used kernel functions for SVM algorithms are linear, polynomial, radial based functions or sigmoid kernels \cite{patle2013svm}, but it is also allowed to design your own kernel as it is easy to make it work with the current libraries that implement SVM algorithms. To select the most appropriate kernel, it is crucial to carefully analyze the dataset beforehand and pre-process it if necessary. Choosing a kernel function is not straightforward, as it significantly affects the results and is highly dependent on the characteristics of the dataset characteristics, such as its linearity or nonlinearity.

The main objective after computing the separation surface  is to maximize the margin between the it and the nearest data points, known as support vectors. The margin maximization is often formulated as an optimization problem using the kernel function defined in Eq. \eqref{Kernel}. Using a generic kernel, the problem can be stated as follows: we want to maximize Eq. \eqref{First_Eq_L_Mul} subjected to Eq. \eqref{Second_Eq_L_Mul}.
\begin{equation} \label{First_Eq_L_Mul}
    \mathcal{F} (\bm{\alpha})
 = \sum_{i = 1}^{n} \alpha_{i} - \frac{1}{n} \sum_{i,j = 1}^{n} \alpha_{i} \alpha_{j} y_{i} y_{j} \mathcal{K}(x_{i}, x_{j}) \;,
\end{equation}
\begin{equation} \label{Second_Eq_L_Mul}
    \sum_{i = 1}^{n} \alpha_{i} y_{i} = 0 \;\;,\;\; 0 \leqslant \alpha_{i} \leqslant\mathcal{C} \;,
\end{equation}
where $\mathcal{C}$ is a constant known as the regularization parameter, which manages the trade-off between achieving high accuracy and ensuring that the model is general enough to be used with new data. The $\alpha_{i}$ in Eq. \eqref{First_Eq_L_Mul} and Eq. \eqref{Second_Eq_L_Mul} are the Lagrange multipliers that are introduced to solve the constrained optimization problem as the conditions for it to be optimum (Karush-Kunh-Tucker conditions) involve them. After solving the problem we arrive at the decision function that allow us to classify a new point $x$ in the featured space:
\begin{equation}
    f(x) = \text{sgn} \left( \sum_{i = 1}^{n} \alpha_{i} y_{i} \mathcal{K}(x,x_{i}) + b \right) \;,
\end{equation}
where $\text{sgn}()$ is the sign function which returns the sign of the expression it has inside once evaluated, and $b$ is called the bias, which is computed using the support vectors ensuring that the solution complies with the margin requirements.

\begin{figure}[!h]
    \centering
    A) \includegraphics[scale = 0.26]{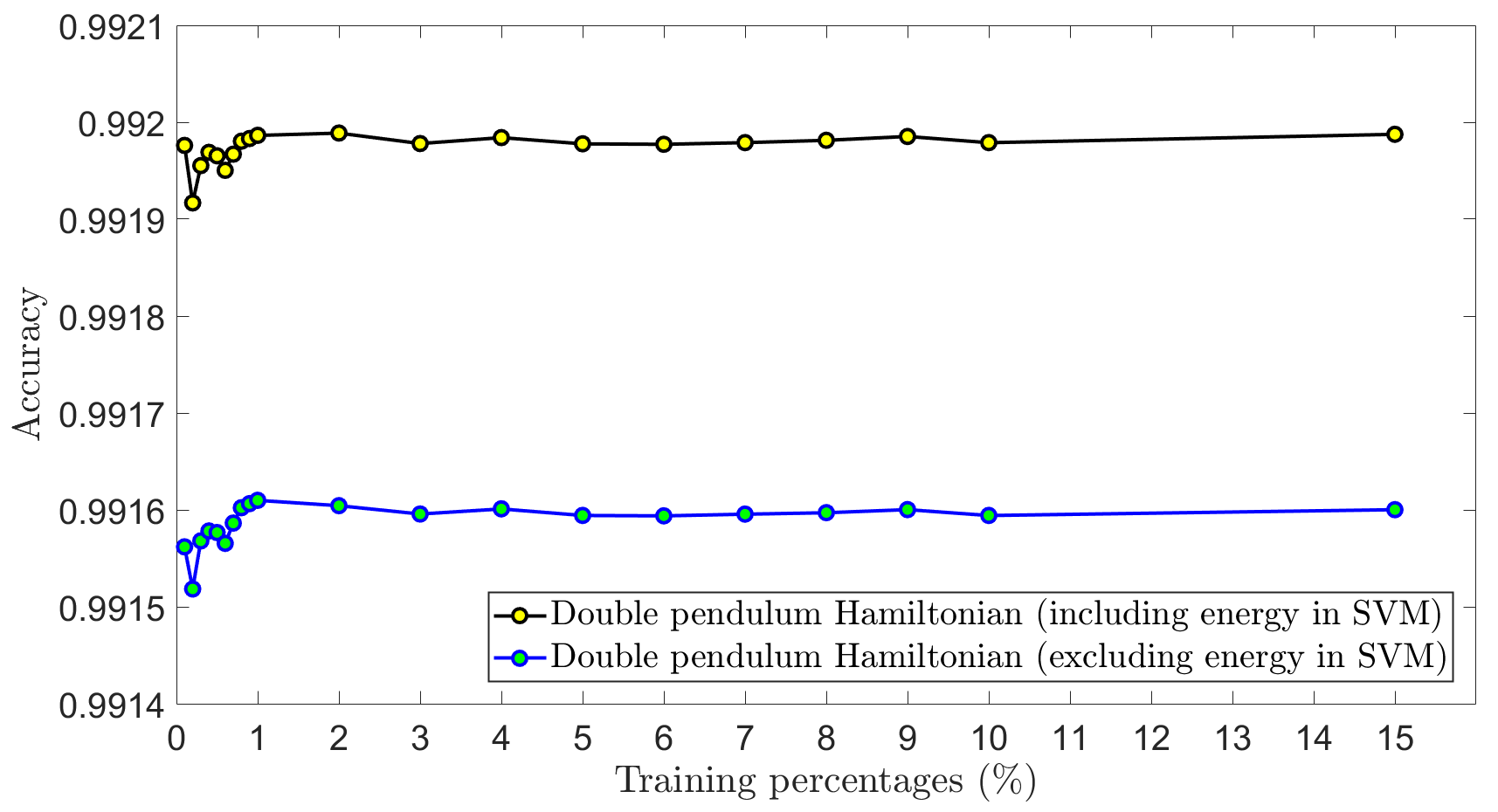} \\
    B) \includegraphics[scale = 0.26]{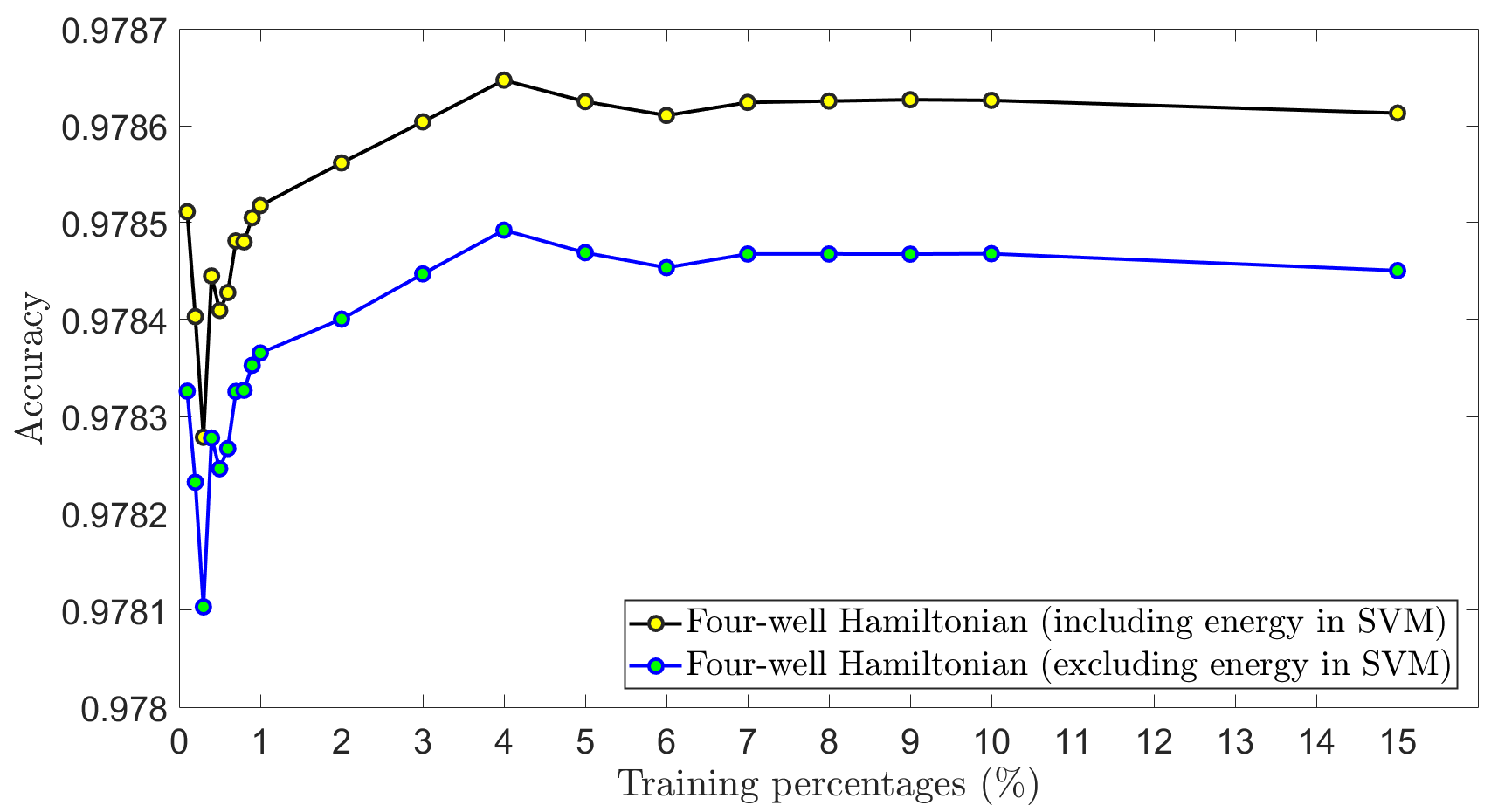}
    \caption{Accuracy of the predictions made by the SVM trained with the $S_L^2(\mathbf{x}_0)$ values for the double pendulum. In A) the accuracy for a different dataset composed of trajectories from the double pendulum is presented while in B) the accuracy obtained for data from the four-well Hamiltonian system is presented.}
	\label{training_svm}
\end{figure}

Support Vector Machines are especially well suited for our case of study due to the fact that our data is binary classified as regular or chaotic. Not only this feature is important, but the main characteristic of our data which allows us to use this simple algorithm is the fact that the value of $S_L^2(\mathbf{x}_0)$ for a regular orbit and for a chaotic one are far way from each other as it can be appreciated in Figs. \ref{Map} and \ref{Henon} or in \cite{jimenez2024pendulum}. This means that we can treat this problem as being linear, and we will have two separated sets of points in the parameter space that we have picked as representative to describe our data. These parameters are the value of the $S_L^2(\mathbf{x}_0)$ indicator and the energy of the system. During the first stages of this project we used both parameters, but after doing some tests we determined that the model which doesn't take the energy into account and only uses the value of $S_L^2(\mathbf{x}_0)$ (or $\log_{10}(S_L^2(\mathbf{x}_0))$) is almost equally efficient, as can be seen in Fig. \ref{training_svm}. We can consider thus the trained SVM model that excludes the energy of the system is much simpler, and hence a better model. For other systems, we have observed that it is more convenient to use the $\log_{10}(S_L^2(\mathbf{x}_0))$ instead of just the value of $S_L^2(\mathbf{x}_0)$ as the logarithmic scale is able to mitigate the effect of outliers in the datasets.

In order to improve computational efficiency, especially when it comes to handle large datasets, parallel computing provides crucial advantages. Utilizing GPUs, particularly those from NVIDIA with CUDA, significantly accelerates the training process compared with the time it takes relying only on CPUs.  Training with CPUs is much slower, but unfortunately finding SVM algorithms which can run in GPU was a challenging task, so we decided to create a neural network that functions equivalently to an SVM using the library PyTorch \cite{NEURIPS2019_9015} and other tools implemented by Scikit-Learn \cite{scikit-learn} to measure the accuracy of our models. We achieved this by using a neural network with one single fully connected layer and a Hinge loss function. This neural network can be formulated as linear transformation where $\boldsymbol{x} \in \mathbb{R}^n$, where $n$ is the number of features, is the input vector and $y \in \mathbb{R}$ is the output, which is computed as:
\begin{equation}
    y = \boldsymbol{w}^{T} \boldsymbol{x} + b \;,
\end{equation}
where $\boldsymbol{w} \in \mathbb{R}^n$ are the weights calculated during the training and $b \in \mathbb{R}$ is the bias vector which is also calculated during the training. After having the output from the layer, we have to apply a prediction method to get the result in the way we want it to be, in our case, binary. To do so, the prediction function we have used is:
\begin{equation}
    \hat{\text{y}} = \text{sgn(y)} \;,
\end{equation}
where the negative values are then mapped to $0$ and the positive ones to the value $1$.

Once the model is defined, we have to train it, for which a stochastic gradient descent algorithm \cite{d8d62392-9a37-31e7-ad3b-37a6f6ee8ef6} was used with $500000$ epochs and the Hinge loss function:
\begin{equation}
    \mathbf{L} = \dfrac{1}{N} \sum_{i = 1}^{N} \text{max} (0, y_{i}(\boldsymbol{w}^{T} \hspace{0.025cm} \boldsymbol{x}_{i} + \text{b}))
\end{equation}
where the number $N$ is the number of samples in each batch. We have selected the Hinge loss function \cite{rojo2018digital} due to its simplicity, efficiency and robustness to noise data even though it has a problem which might be relevant in some cases, and its the fact that it is not differentiable at $0$, which can make it hard to optimize with some of the gradient-based methods.

In order to train our neural network, we used values of the $S_L^2(\mathbf{x}_0)$ indicator calculated for the double pendulum with the goal of classifying later the trajectories for the four-well Hamiltonian, which were initially classified using the algorithm described in \cite{jimenez2024pendulum} to obtain the threshold value for indicators based on Lagrangian Descriptors. The reasons for using this indicator in particular, considering that there are a lot of them, are described in \cite{zimper23}, but it is mainly for it high precision when compared to the SALI method. On the other hand, for the H\'enon-Heiles system and the Chirikov Standard Map the $\log_{10}(S_L^2(\mathbf{x}_0))$ provided much better results during the training and validation processes, as in logarithmic scale some of the outliers of our dataset were no longer a problem, see Tables \ref{svm_vs_sali_results_HH} and \ref{svm_vs_sali_results_StdMap}.

After training the model with the data from the pendulum, we have validated it with the systems already mentioned. For the four-well Hamiltonian, we did this using the value of the indicator $S^2_L$ computed for more that $22$ million trajectories for both model approaches, including and excluding the energy. For the other systems, we computed both indicators, SALI and $S^2_L$ so we can compare with a grown truth method and be sure about how well our models have performed for different systems that only share the feature of having $2$ Degrees of Freedom.

\section{Results} \label{Results}

In this section we will describe the results of the conducted research. We have trained different models using the double pendulum data we simulated for one of our previous papers \cite{jimenez2024pendulum} varying the amount of data we gave the model to be trained with. Two different kinds of models were explored: one was provided with data for the chaos indicator $S_L^2(\mathbf{x}_0)$, computed for an integration time of $700$ time units, and the energy of the system and the other was only provided with the values of the chaos indicator. In Fig.\ref{training_svm} A) the accuracy of both models for different amount of training data is showed, were it is easily appreciated that both models perform incredibly well when it comes to classifying the double pendulum data they were not trained with. This fact is telling us that the model is able to learn the system's behavior with few data from it. Furthermore, both approaches work incredible well, so it is logical to choose the model which only uses the values of the chaos indicator as it is simpler and because it will allow us to use it also for discrete maps, in which it is not possible to define an energy itself.

In order to test if our models, which have only been trained with the data from one Hamiltonian system, were able to predict the behavior of other system's trajectories we used a four-well potential Hamiltonian \cite{garcia2020exploring} integrated for $700$ units of time to obtain the chaos indicators derived from Lagrangian descriptors. In Fig.\ref{training_svm} B), the accuracy for the two types of models is presented. In this graph we observe that the accuracy is also very high, as it was for the double pendulum, but it is important to note that the models have never seen data coming from this system, but they were able to predict the behavior of almost all of the trajectories, which in this case were 22 million. The classification of these trajectories was originally done with the algorithm developed in \cite{jimenez2024pendulum}, which has proven to be consistent and accurate. Again, as we saw for the double pendulum, the model trained with only the value of the chaos indicator performs a little bit worse, even though both are very close for all the training percentages.

\begin{figure}[!h]
    \centering
    A) \includegraphics[scale=0.245]{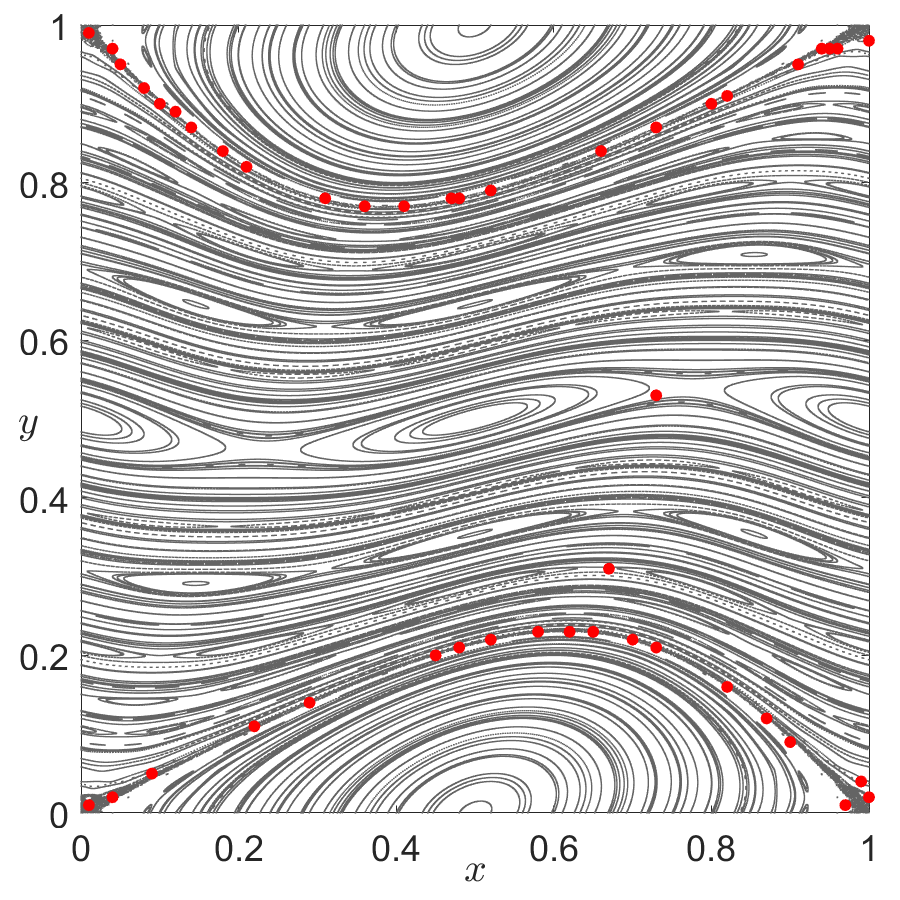}
    B) \includegraphics[scale=0.245]{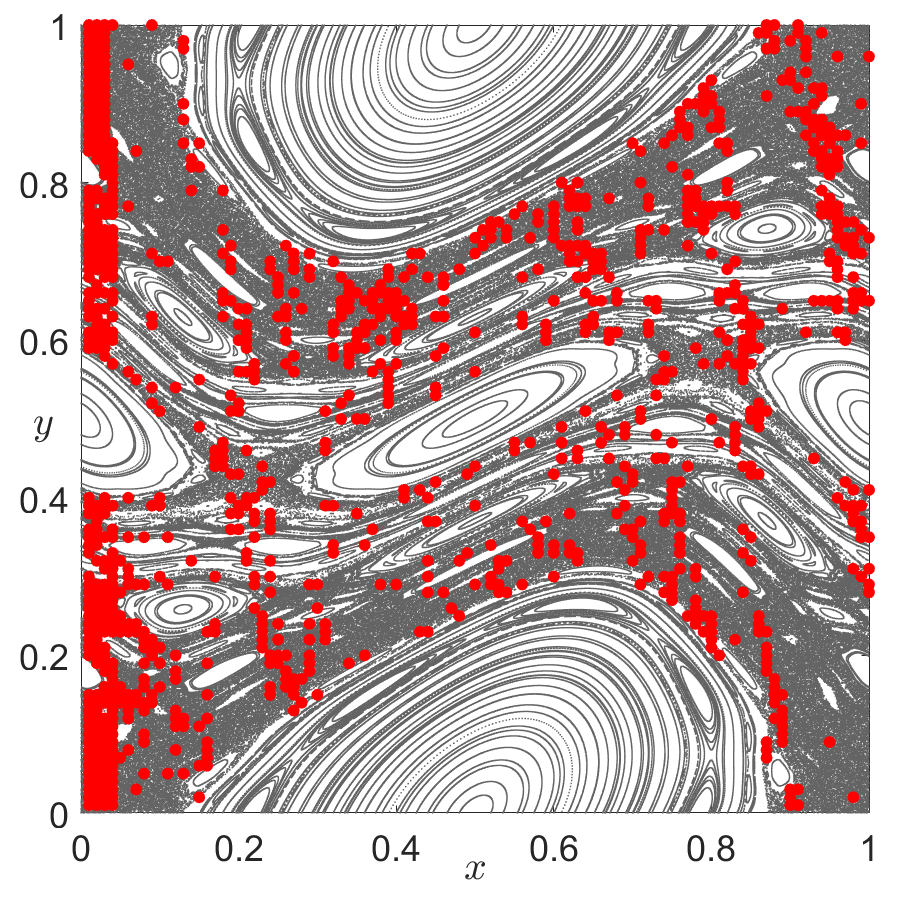} \\
    C) \includegraphics[scale=0.245]{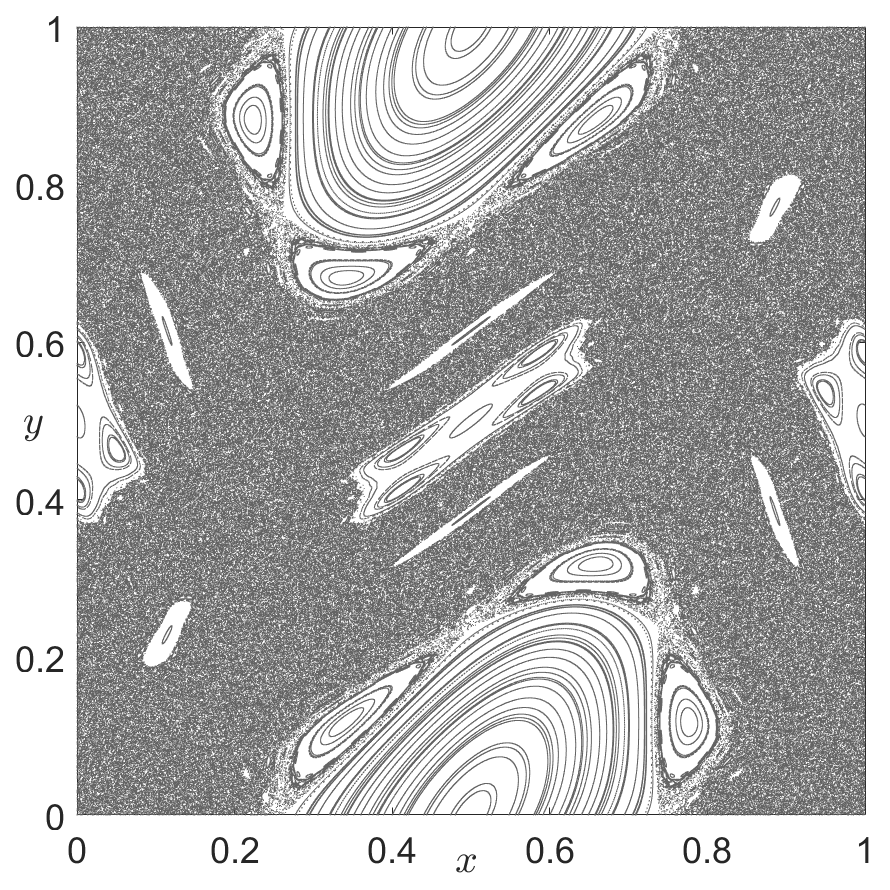}
    \caption{Poincar\'e maps and misclassified initial conditions (red dots) by the SVM model for the Standard Map. A) $K=0.5$; B) $K=0.971635$; C) $K=1.5$. The SVM model was trained with the $\log_{10}(S_L^2)$ dataset for the double pendulum Hamiltonian.}
    \label{fig:miss_std_map}
\end{figure}

\begin{figure}[!h]
    \centering
    A) \includegraphics[scale=0.245]{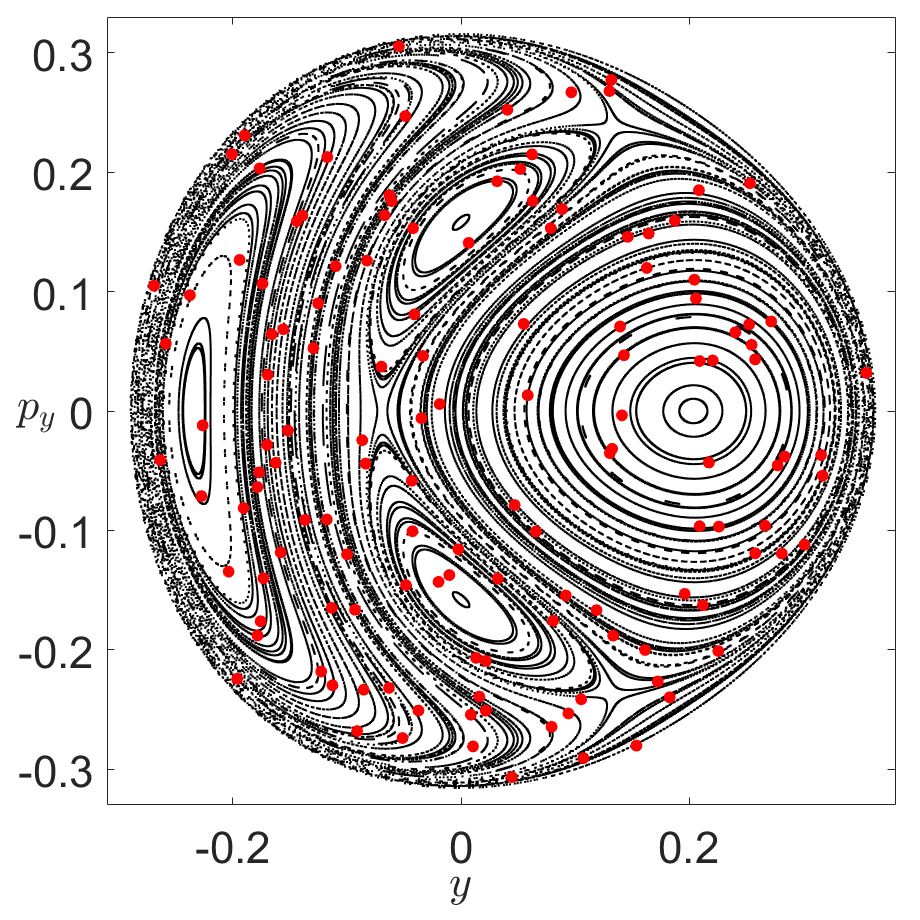}
    B) \includegraphics[scale=0.245]{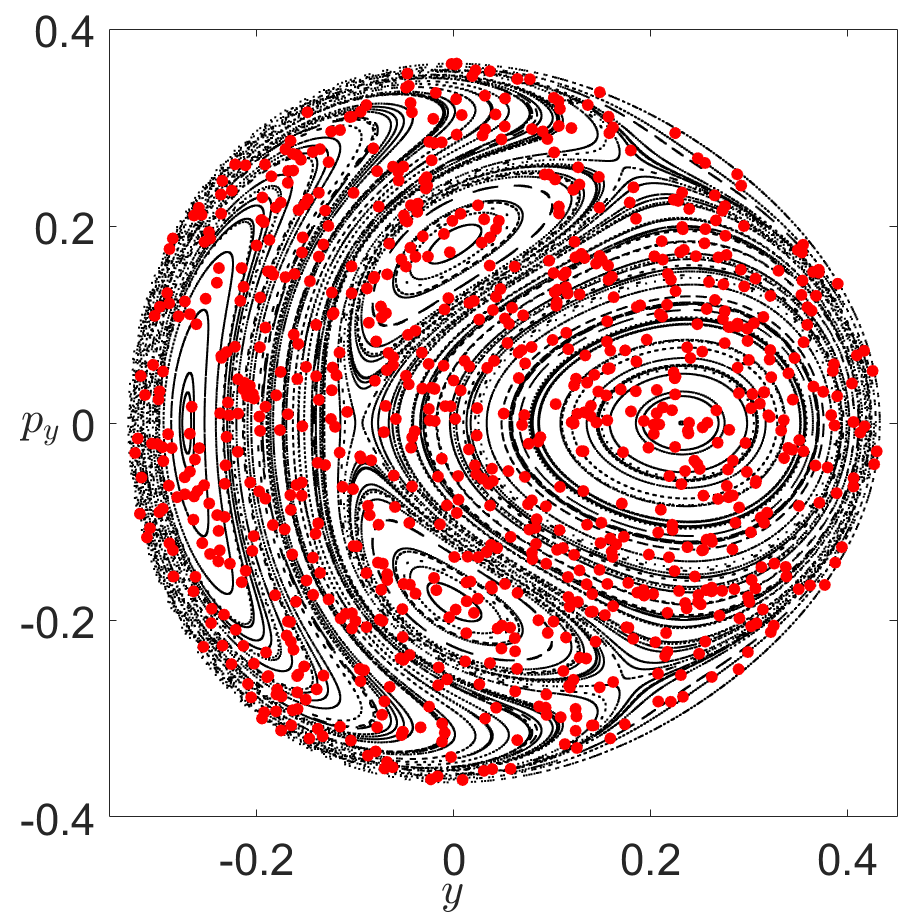} \\
    C) \includegraphics[scale=0.245]{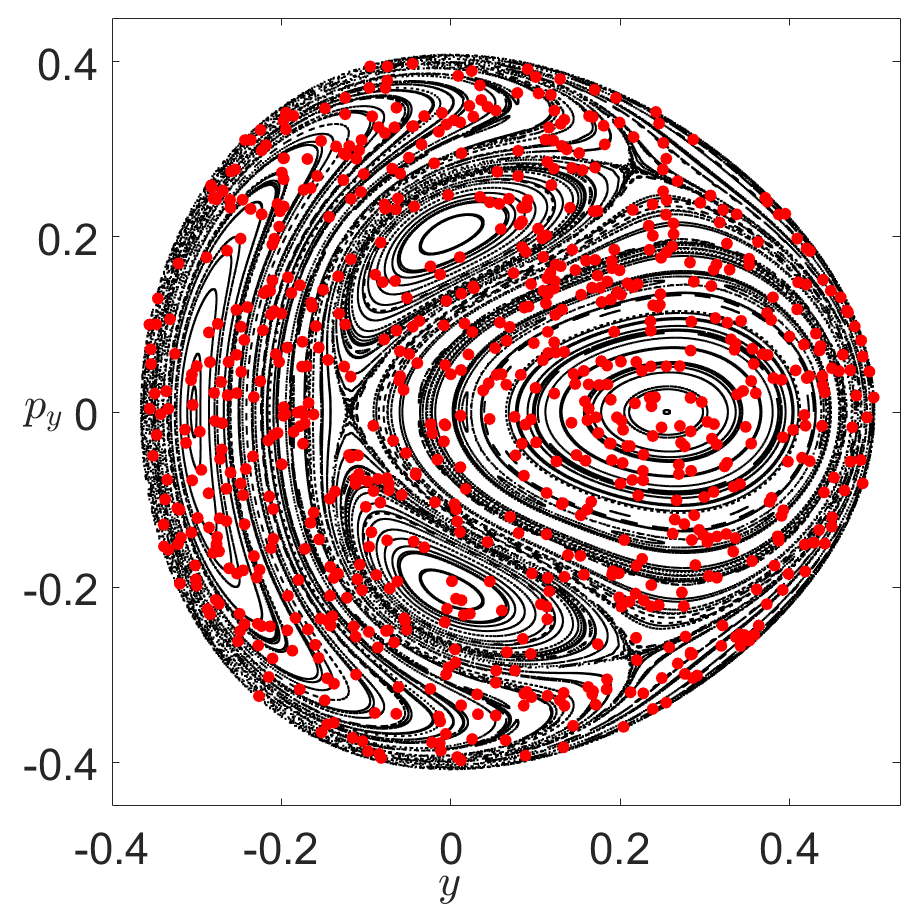}
    D) \includegraphics[scale=0.245]{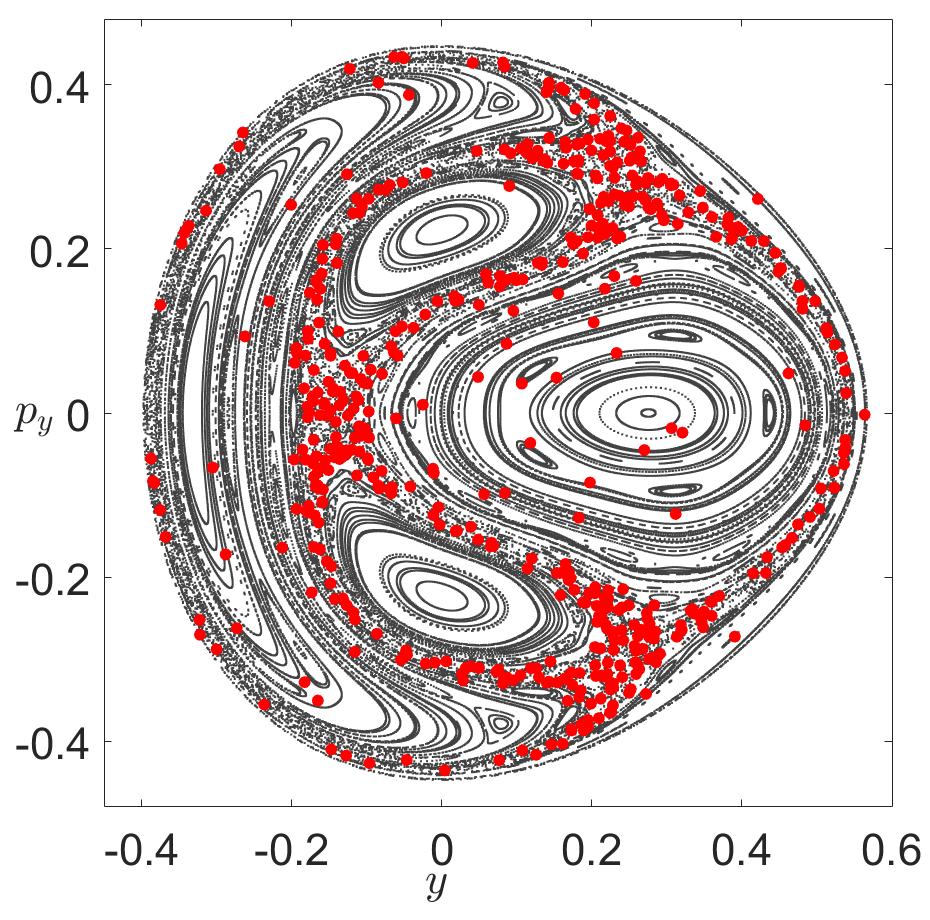} \\
    E) \includegraphics[scale=0.245]{missClass_HH_1div8} 
    \caption{Poincar\'e maps and misclassified initial conditions (red dots) by the SVM model for the H\'{e}non-Heiles Hamiltonian. A) $\mathcal{H}=1/20$; B) $\mathcal{H}=1/15$; C) $\mathcal{H}=1/12$; D) $\mathcal{H}=1/10$; E) $\mathcal{H}=1/8$. The SVM model was trained with the $\log_{10}(S_L^2)$ dataset for the double pendulum Hamiltonian.}
    \label{fig:miss_HH}
\end{figure}

\begin{figure}[!h]
    \centering
    A) \includegraphics[scale=0.245]{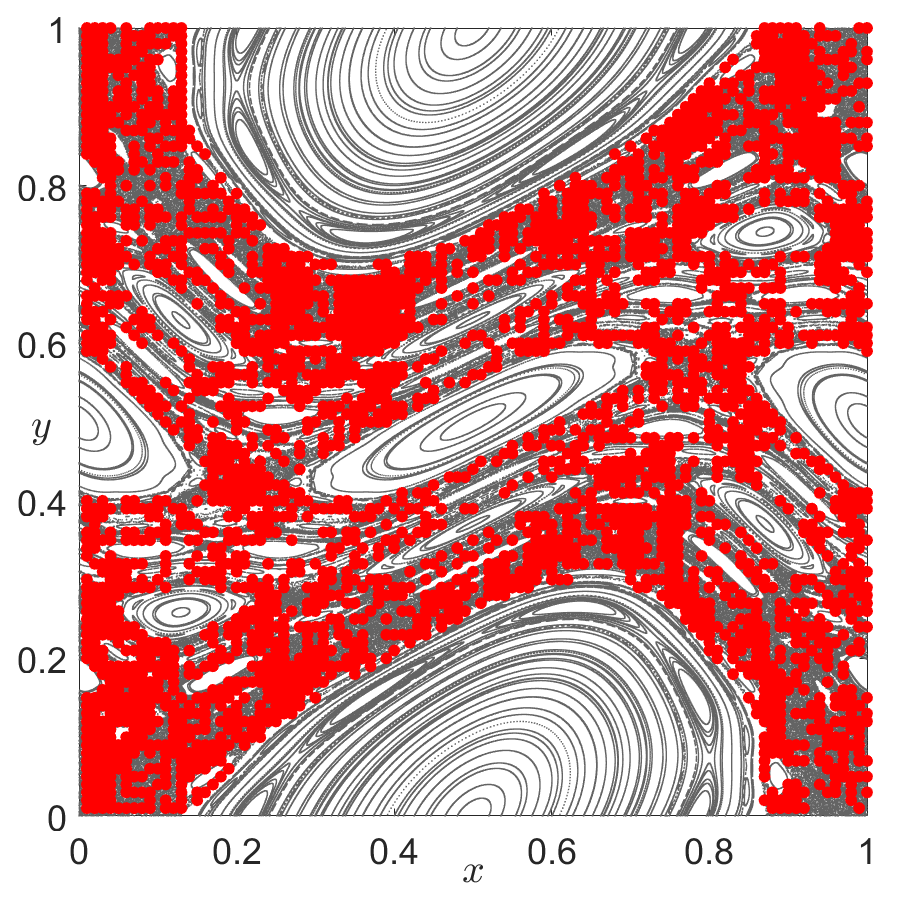}
    B) \includegraphics[scale=0.245]{missClass_StdMap_K_0p971635}
    \caption{Poincar\'e maps and misclassified initial conditions (red dots) by the SVM model for the Standard Map with $K=0.971635$. A) The SVM model was trained with the $S_L^2$ dataset for the double pendulum Hamiltonian; B) The training of the SVM was carried out with the $\log_{10}(S_L^2)$ values.}
    \label{fig:miss_std_map_S_and_logS}
\end{figure}

\begin{figure}[!h]
    \centering
    A) \includegraphics[scale=0.245]{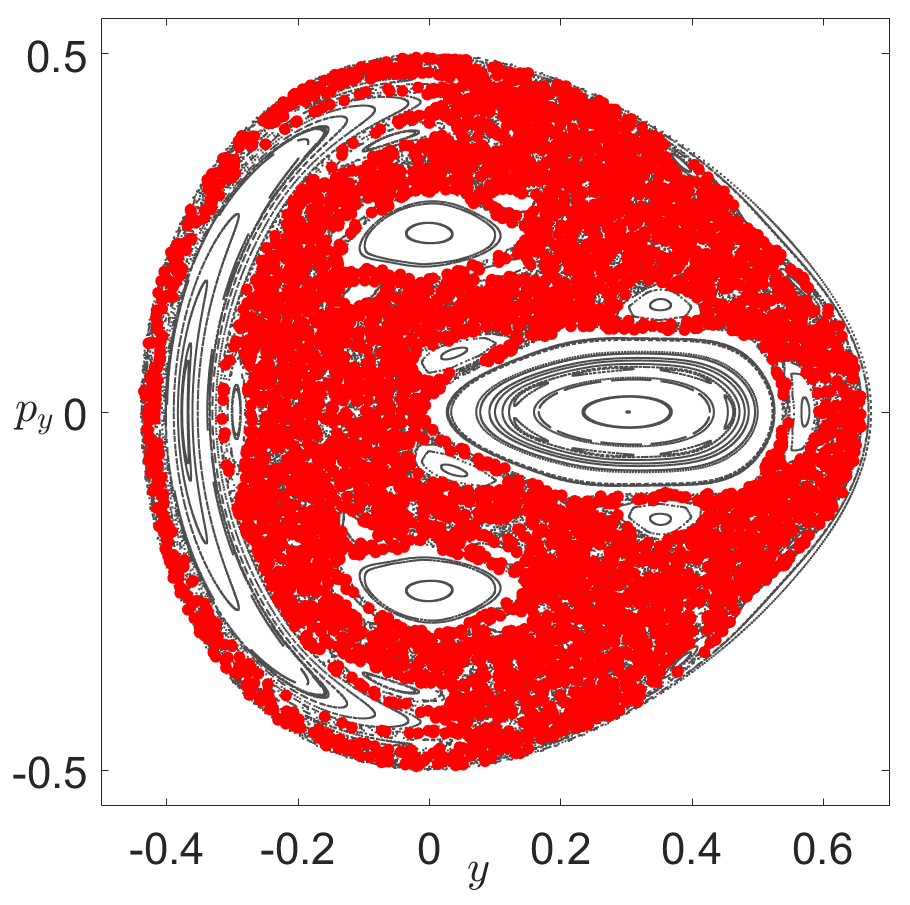}
    B) \includegraphics[scale=0.245]{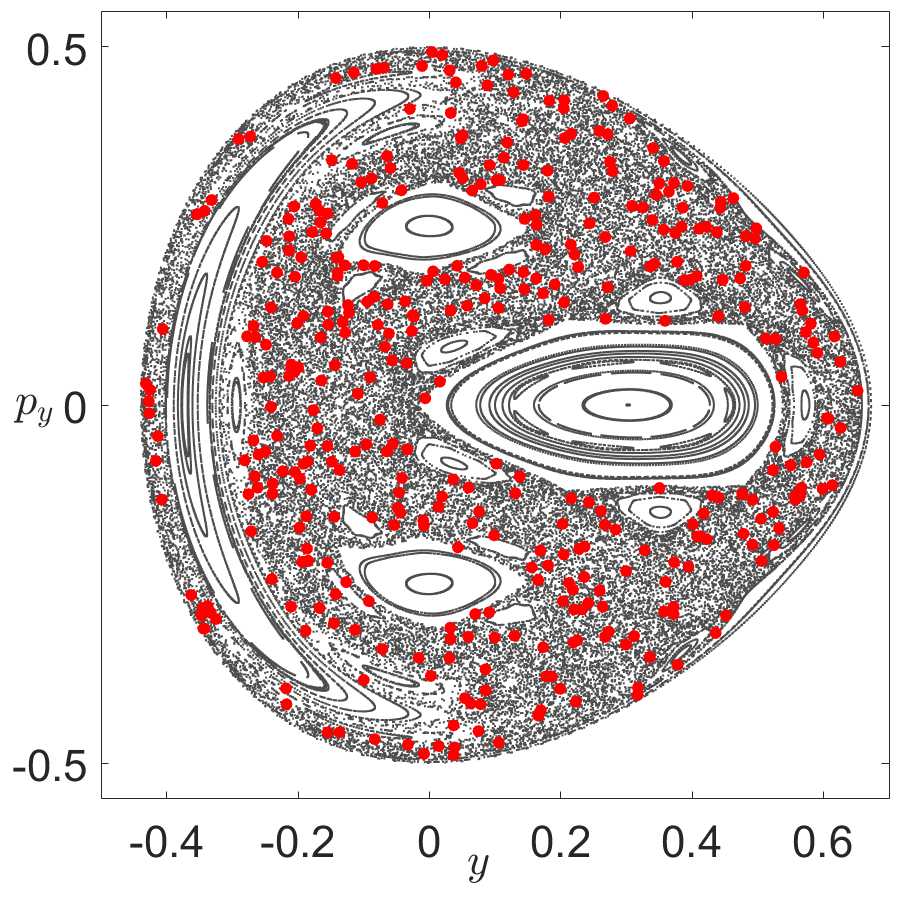}
    \caption{Poincar\'e maps and misclassified initial conditions (red dots) by the SVM model for the H\'enon-Heiles Hamiltonian with energy $\mathcal{H} = 1/8$. A) The SVM model was trained with the $S_L^2$ dataset for the double pendulum Hamiltonian; B) The training of the SVM was carried out with the $\log_{10}(S_L^2)$ values.}
    \label{fig:miss_HH_S_and_logS}
\end{figure}

To show how powerful these models are despite of their simplicity, we have trained a model with just the $10 \%$ of the data from the double pendulum but using the logarithm of the indicator instead of the value itself due to the reasons explained in section \ref{Methodology}, and predicted the behavior of trajectories in the Hénon-Heiles system \cite{henon1964applicability} with $50000$ initial conditions, $10000$ for each of the $5$ energies simulated and with $30000$ for the Chirikov Standard Map \cite{Chirikov1979}, with another $10000$ for each value of the parameter $K$. For both of these systems, the "manual" classification was performed using the SALI method as shown in Fig. \ref{Henon} B) and Fig. \ref{Map} B), which we consider to be the ground truth chaos indicator, and then compared with the predictions of the model. In Tables \ref{svm_vs_sali_results_HH} and \ref{svm_vs_sali_results_StdMap} we show the accuracy that the SVM model has when compared with the SALI indicator to classify the chaotic or regular nature of trajectories. As one can clearly see from these results, it is very high for all of the presented cases, so we can conclude that the model is performing in an excellent way. Note that the SVM model trained with the values of $\log_{10}(S^2_L)$ performs better than that using the $S_L^2$ indicator. In Figs. \ref{fig:miss_std_map} and \ref{fig:miss_HH} we depict with red dots, overlaid with the corresponding Poincar\'e map, the initial conditions that have been misclassified by the SVM model (when compared with the ground truth classification provided by the SALI indicator), both for the Standard Map and for the H\'enon-Heiles Hamiltonian. In order to compute the SALI indicator, we have to integrate the variational equations of the system during enough time so that the chaotic and regular behaviors will be correctly identified. In this case, we used an integration time for SALI of $10^5$ time units, while for the chaos indicators derived for Lagrangian descriptors, explained in section \ref{LD_ind}, we used a time of $10^3$ for the H\'enon-Heiles system and $5\times10^3$ units of time for the Chirikov Standard Map. These misclassified initial conditions are usually lying in the chaotic region of the phase space for both systems, which was somehow expected as the chaotic behavior is always harder to analyze. But sometimes this is not correct and the misclassified initial conditions lay all over the Poincaré section as it's shown in Fig. \ref{fig:miss_std_map_S_and_logS} and Fig. \ref{fig:miss_HH_S_and_logS} without apparent reason. Exploring this more deeply, we can notice that when we use the value of the indicator $S^2_L$ instead of $\log_{10}(S^2_L)$, most of the chaotic trajectories are misclassified. 

Our results show that a simple model is perfectly suitable for the purpose of detecting chaos and regularity in Hamiltonian systems, and it is able to do it with astonishing results compared to more complex techniques that have been applied for the same problem \cite{Barrio2023, Seok2020}, proving that SVMs combined with Lagrangian descriptors are a very powerful tool in the field of non-linear dynamics and chaotic dynamics.

\begin{table}[!h]
    \centering
    \begin{tabular}{| c | c | c |}
    \hline
    \multicolumn{3}{|c|}{H\'enon-Heiles Hamiltonian} \\ \hline
    $\mathcal{H}$ & \makecell{Accuracy SVM vs. SALI \\ (trained with $S_L^2$)} & \makecell{Accuracy SVM vs. SALI \\ (trained with $\log_{10}(S_L^2)$)}\\ \hline 
    $1/20$ & 100\% & 98.66\% \\ \hline
    $1/15$ & 100\% & 91.3\%\\ \hline
    $1/12$ & 100\% & 91\% \\ \hline
    $1/10$ & 86.38\% & 94.67\% \\ \hline
    $1/8$ & 40.22\% & 96.1\% \\ \hline
    \end{tabular}
    \caption{Accuracy of the predictions made by the SVM model for the H\'enon-Heiles Hamiltonian, compared to the ground truth classification provided by the SALI indicator. The model was trained with the double pendulum dataset, and the energy was not included as a parameter during the training process.} 
    \label{svm_vs_sali_results_HH}
\end{table}

\begin{table}[!h]
    \centering
    \begin{tabular}{| c | c | c |}
    \hline
    \multicolumn{3}{|c|}{Chiirokov Standard Map} \\ \hline
    $K$ & \makecell{Accuracy SVM vs. SALI \\ (trained with $S_L^2$)} & \makecell{Accuracy SVM vs. SALI \\ (trained with $\log_{10}(S_L^2)$)}\\ \hline 
    0.5 & 98.52\% & 99.55\% \\ \hline
    0.971635 & 68.16\% & 89.55\% \\ \hline
    1.5 & 96.1\% & 100\% \\ \hline
    \end{tabular}
    \caption{Accuracy of the predictions made by the SVM model for the Chirikov Standard Map, compared to the ground truth classification provided by the SALI indicator. The model was trained with the double pendulum dataset, and the energy was not included as a parameter during the training process.} 
    \label{svm_vs_sali_results_StdMap}
\end{table}


\section{Conclusions} \label{Conclusions}

In this paper we have investigated how a very simple model, a binary classifier implemented as a SVM, is able to predict the regular and chaotic behavior of trajectories in Hamiltonian systems with great performance. Our findings have revealed that one Hamiltonian system, in our particular case the double pendulum, encodes enough information to predict and understand the behavior of other Hamiltonian systems which are not related to the original one.

Our first approach for the construction of a SVM model was to provide data for the energy of the system and the value of the $S_L^2(\mathbf{x}_0)$, obtaining excellent results. Shortly after, we decided to try and get rid off the energy as something we need to classify a trajectory, providing a simpler model. The accuracy of these models, despite being slightly lower, is still excellent as we have already demonstrated when testing both approaches in the four-well Hamiltonian system. After, we trained another type of SVM model in which the input data was the $\log_{10}(S_L^2(\mathbf{x}_0))$. This last option has been proven to be more accurate when tested in the H\'enon-Heiles system and the Chirikov Standard Map due to the presence of outliers in the data, which in logarithmic scale were no longer a problem and allowed us to obtain excellent results for the classification of trajectories.

Perhaps, one of the most interesting findings of this work is that, due to the way Lagrangian descriptors are constructed and how the chaos indicators derived from them are defined and computed for both, continuous and discrete systems, we were able to classify the dynamics of multiple initial conditions in a discrete system using a model that was only trained with data from a continuous system.

Our novel approach to chaos classification using SVMs, neural networks and Lagrangian descriptors gives significant improvement to the field as the previous ways to classify initial conditions, which is the same as obtaining the value of the corresponding chaos indicator that separates regular conditions from chaotic ones, were mainly two: doing it by hand or using complex and slow algorithms. In both cases, we needed to analyze data sampled from a large range of energies, where for each energy we used $10^4$ initial conditions, making this a computationally demanding and slow process. With our SVM model it is enough to simulate the time evolution of the trajectory associated to the initial condition whose regular or chaotic nature we aim to distinguish, and this allows us to simultaneously apply this procedure to classify an arbitrary number of initial conditions. This reduces considerably the time required to investigate regularity and chaos of Hamiltonian systems.

The development and implementation of machine learning models aimed at the classification of the regular and chaotic nature of trajectories, trained with data obtained from Lagrangian descriptors, opens up many interesting avenues for future research. The use of LDs for this task has the important advantage that, in comparison to other chaos indicators, their computation only requires the integration of trajectories using the equations of motion. This work constitutes a first step towards achieving this goal, and shows how a simple approach as SVMs are, when compared to other Deep and Machine Learning techniques such as Convolutional Neural Networks, can be easily applied to this problem in order to  provide accurate results compared to the more complex algorithms that are used in the literature for distinguishing between regular and chaotic behavior. 

Our future research will focus on developing more complex models to assess if accuracy can be improved, maintained, or reduced. We will also use more data from different systems to create more general models, experiment with different training algorithms, and test other chaos indicators derived from Lagrangian descriptors. Additionally, we plan to conduct a more exhaustive study of the different parameters that can be considered for an SVM, such as the value of $\mathcal{C}$ and the number of epochs used during the training stage.


\section*{Data availability}

The data that support the findings of this study are available from the corresponding author upon reasonable request.

\section*{CRediT authorship contribution statement}

\textbf{Javier Jim\'enez L\'opez:} Conceptualization, Methodology, Data curation, Formal analysis, Funding acquisition, Investigation, Resources, Software, Validation, Visualization, Writing - original draft, Writing - review \& editing.

\noindent
\textbf{V\'{i}ctor J. Garc\'{i}a-Garrido:} Conceptualization, Methodology, Formal analysis, Funding acquisition, Investigation, Project administration, Software, Supervision, Validation, Visualization, Writing - original draft, Writing - review \& editing.

\begin{appendices}
\label{Appendix}

\section{Model Hamiltonian Systems} 

This appendix briefly describes the mathematical formulation of the Hamiltonian systems we have used for this work. For our analysis we have considered a two-dimensional symplectic map, the well known Chirikov Standard Map (CSM), and three different Hamiltonian systems with 2 degrees of freedom (DoF).

The Chirikov Standard Map \cite{Chirikov1979}, also known as the Standard Map, was introduced by the Russian physicist Boris Chirikov in 1969 for understanding the transition from regular (integrable) to chaotic (non-integrable) motion in Hamiltonian systems. This model system was developed in the context of plasma physics and accelerator dynamics, particularly in relation to the behavior of particles in magnetic fields and the stability of motion in nonlinear dynamical systems. It has since become a fundamental example in the study of chaotic dynamics, and has been widely used to illustrate the principles of chaos theory in a variety of physical and mathematical contexts \cite{Meiss2008}. The Standard Map is a two-dimensional area-preserving map \cite{Meiss1992} which displays a rich variety of dynamical behaviors. Given a general planar discrete dynamical system defined by the evolution rule:
\begin{equation}
    \mathbf{z}_{n+1} = \mathbf{f}(\mathbf{z}_n) \;\;,\;\; n \in \mathbb{N} \cup \lbrace 0 \rbrace \;,
\end{equation}
with $\mathbf{z}_n = (x_n,y_n) \in \mathbb{R}^2$, the CSM is given by the system of difference equations:
\begin{equation}
    \begin{cases}
        x_{n+1} = x_n + y_{n+1} \\[.1cm]
        y_{n+1} = y_n + \dfrac{K}{2\pi} \sin(2\pi x_n)
    \end{cases} (\text{mod } 1) \;,
    \label{std_map}
\end{equation}
where $x_n \in [0,1]$ and $y_n \in [0,1]$ represent, respectively, the position and momentum of the particle at a discrete time step $n$, and $K$ is the parameter that controls the strength of the nonlinearity and the degree of chaos in the system. In our study, we have selected the values of $K = 0.5$, $K = 0.971635$ and $K = 1.5$, and for each of these cases we have randomly sampled an ensemble of $10^4$ initial conditions on the square $[0,1]\times[0,1]$.

\begin{figure}[!h]
    \centering  
    A) \includegraphics[scale = 0.26]{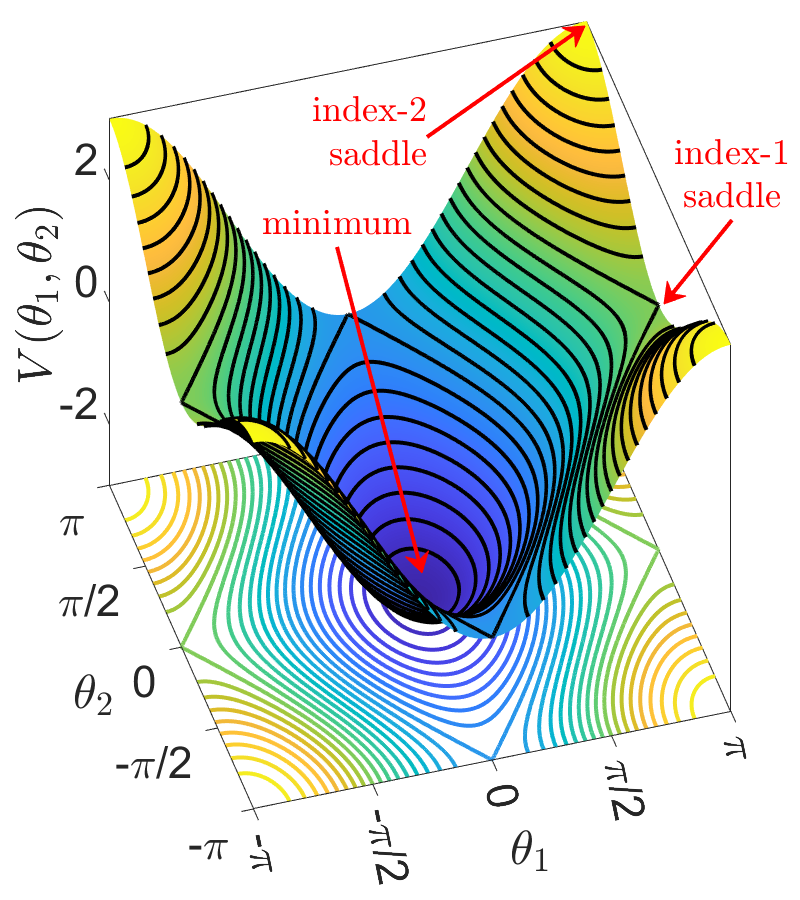}
    B) \includegraphics[scale = 0.25]{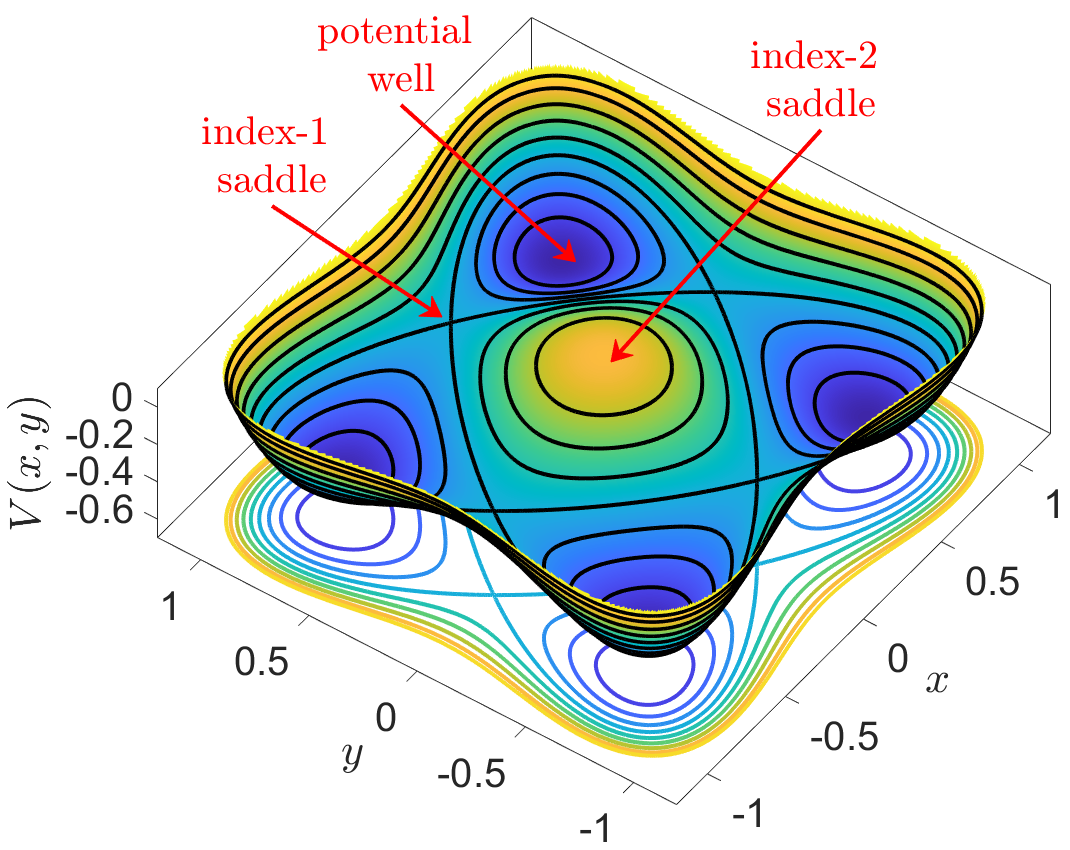} \\
    C) \includegraphics[scale = 0.27]{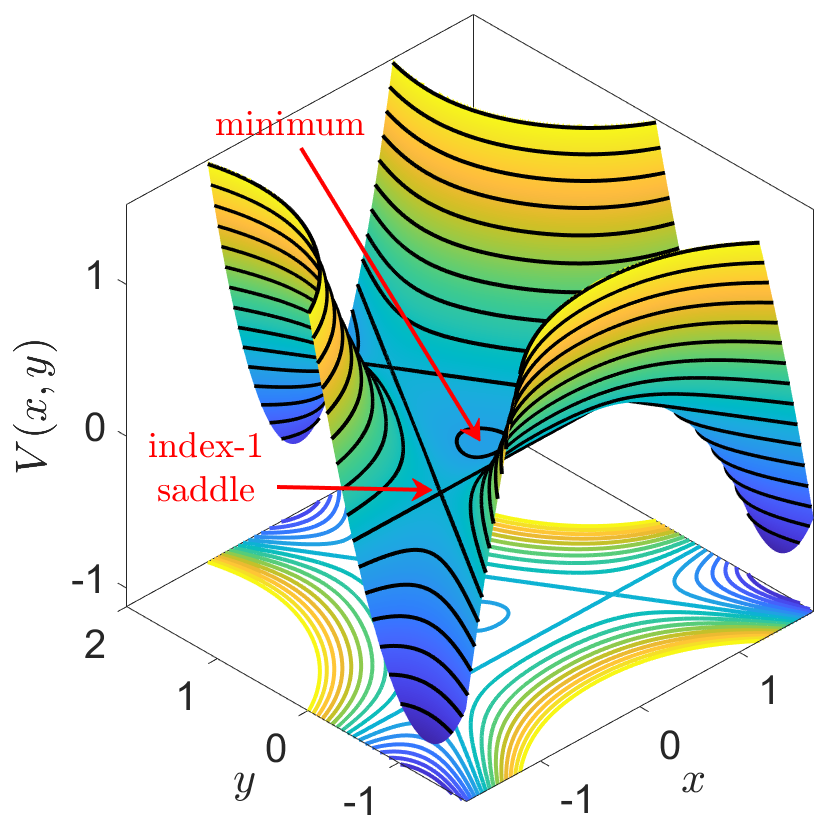}
    \caption{Potential energy function and equipotential contours for the three Hamiltonian systems with two degrees of freedom that we have used in ths paper. A) PES in Eq. \eqref{pes_dpend} for the double pendulum system; B) PES in Eq. \eqref{pes_fourwells} for the four-well system; C) PES in Eq. \eqref{pes_HenonHeiles} for the H\'enon-Heiles system.}
    \label{PESs}
\end{figure}

The Hamiltonian system with 2 DoF that we have used to train and validate our SVM model is the well-known double pendulum \cite{korsch2007chaos}, which is a mechanical system consisting of two pendulums attached end to end. If we consider that the masses and lengths of the pendulums are, respectively, $m_1$, $m_2$ and $l_1$, $l_2$, then one can define the dimensionless parameters, $\alpha = l_1/l_2$ and $\sigma = m_1/m_2$, and study the chaotic and regular dynamics of this system in terms of these model parameters. In fact, this task was recently undertaken in \cite{jimenez2024pendulum}, where the classification trajectories according to their chaotic and regular nature was carried out by means of chaos indicators based on Lagrangian descriptors. In \cite{jimenez2024pendulum} it was shown that the Hamiltonian for the double pendulum can be written compactly in the form:
\begin{equation} 
\mathcal{H}(\boldsymbol{\theta},\mathbf{p}) = \dfrac{1}{2} \mathbf{p}^T M^{-1}(\cos \Delta \theta) \, \mathbf{p} - \alpha(1+\sigma) \cos \theta_1 - \cos \theta_2 \;,
\label{Ham_dpend}
\end{equation}
where $\boldsymbol{\theta} = (\theta_1,\theta_2)$ are the configuration space coordinates representing the angle that each of the pendulums make with the vertical, $\Delta \theta = \theta_2 - \theta_1$ is the angular difference, and $\mathbf{p} = (p_1,p_2)$ is the conjugate momenta vector. The inverse of the mass matrix that appears in the kinetic energy term of the Hamiltonian function is given by:
\begin{equation}
M^{-1}(x) = \dfrac{1}{1+\sigma - x^2}\begin{bmatrix}
	\dfrac{1}{\alpha^2} & -\dfrac{x}{\alpha} \\[.35cm]
	-\dfrac{x}{\alpha} & 1+\sigma
\end{bmatrix} \;.
\end{equation}
Moreover, the potential energy surface (PES) is defined by:
\begin{equation}
    V(\theta_1,\theta_2) = - \alpha(1+\sigma) \cos \theta_1 - \cos \theta_2 \;.
    \label{pes_dpend}
\end{equation}
and, as we illustrate in Fig. \ref{PESs} A), it has a local minimum at the origin (with center$\times$center stability), two index-1 saddles (with saddle$\times$center stability) and an index-2 saddle (with saddle$\times$saddle stability). Hamilton's equations for this system can be written as follows:
\begin{equation}
\begin{cases}
    \dot{\boldsymbol{\theta}} = M^{-1}(\cos \Delta \theta) \, \mathbf{p} \\ 
    \dot{\mathbf{p}} = \dfrac{\sin \Delta \theta}{2} \mathbf{p}^T C(\cos \Delta \theta) \, \mathbf{p} \begin{bmatrix}
		1 \\
		-1
	\end{bmatrix} - \begin{bmatrix}
	\alpha(1+\sigma) \sin\theta_1 \\[.1cm]
	\sin\theta_2
	\end{bmatrix}
\end{cases} .
\label{ham_eqs}
\end{equation}
with $C(x) = -M^{-1}(x) \dfrac{dM}{dx} M^{-1}(x)$. The dataset that we have used to develop the SVM model (training plus validation) is the one we obtained in \cite{jimenez2024pendulum} with the goal of analyzing the chaotic fraction of the phase space in the double pendulum. This dataset has been generated in the following way. We have selected the parameters values $\alpha_i = 2^{i}$ and $\sigma_j=2^{j}$, where $i,j \in \left\lbrace -4,-3,-2,-1,0,1,2,3,4 \right\rbrace$, and, in each of these cases, we have considered $170$ different energy levels ($40$ of them distributed uniformly from the energy of the local minimum up to the energy of the index-2 saddle, and the other $130$ taken above the energy of the index-2 saddle with a unit step). For each numerical experiment we have calculated the chaos indicators based on LDs for a random ensemble of $10^4$ initial conditions, using an integration time of $\tau = 700$. For more details, see \cite{jimenez2024pendulum}.

Another Hamiltonian system that we have used to validate our SVM model was investigated in \cite{garcia2020exploring}, and is inspired from studies on double proton transfer chemical reactions \cite{Smederchina2008}. The Hamiltonian function is:
\begin{equation}
    \mathcal{H}(x,y,p_{x},p_{y}) = \dfrac{p_{x}^2}{2} + \dfrac{p_{y}^2}{2} + V(x,y) \;,
    \label{fourwells_ham}
\end{equation}
where $\mathcal{V}(x,y)$ is the potential energy surface of the model, given by:
\begin{equation}
    V(x,y) = x^4 - \alpha x^2 - \delta x + y^4 - y^2 + \beta x^2 y^2 \;,
    \label{pes_fourwells}
\end{equation}
and $\alpha,\beta,\delta > 0$ are the model parameters. The potential energy landscape of this system is mainly characterized by four wells (center$\times$center stability), four index-1 saddles (saddle$\times$center stability) and an index-2 saddle (saddle$\times$saddle stability), see Fig. \ref{PESs} B). The parameters control the depth of the potential wells ($\alpha$), the coupling of the $x$ and $y$ DoF ($\beta$) and the asymmetry of the potential in the $x$-direction ($\delta$). It is important to remark that by varying the value of $\delta$, the geometrical features that comprise the topography of the PES can undergo drastic changes, as shown in \cite{garcia2020exploring}, and this process has a major impact on the phase space structures of the system. Since our Hamiltonian has 2 DoF and energy is conserved, dynamics takes place on a three-dimensional energy hypersurface of the four-dimensional phase space. Hamilton's equations of motion for this system have the form.
\begin{equation}
\begin{cases}
\dot{x} = \dfrac{\partial H}{\partial p_x} = p_x \\[.4cm]
\dot{y} = \dfrac{\partial H}{\partial p_y} = p_y \\[.4cm]
\dot{p}_x = -\dfrac{\partial \mathcal{H}}{\partial x} = -4x^3 + 2\alpha x + \delta - 2\beta y^2 \\[.4cm]
\dot{p}_y = -\dfrac{\partial \mathcal{H}}{\partial y} = -4y^3 + 2y - 2\beta x^2
\end{cases}
\;.
\end{equation}
The dependence of the potential energy in Eq. \eqref{pes_fourwells} on three parameters has allowed us to easily change the geometrical features of the PES, and thus we have been able to generate a wide range of scenarios with different distributions of chaos and regularity throughout the underlying phase space of this model. It is interesting to note here that, as the model parameter $\delta$ is varied, the geometry of the PES can change drastically, giving rise to different types of bifurcations in the phase space. In our simulations for this system we have used $300$ energy levels measured from that of the minimum of the potential, and we have considered several cases for the parameter values, which we summarized in the table below.
\begin{table}[!h]
    \centering
    \begin{tabular}{| c | c | c | c |}
    \hline
     & $\alpha$ & $\beta$ & $\delta$ \\ \hline 
    Case 1 & 0.75 & 0.25 & 0.5 \\ \hline 
    Case 2 & 1 & 0.25 & 0.1 \\ \hline
    Case 3 & 1 & 1 & 0.1 \\ \hline
    Case 4 & 2 & 0.1 & 0.1 \\ \hline
    Case 5 & 0.5 & 0.75 & 0 \\ \hline
    Case 6 & 1 & 2 & 0 \\ \hline
    Case 7 & 0.01 & 0.01 & 0 \\ \hline
    Case 8 & 0.01 & 0.01 & 0.75 \\ \hline
    \end{tabular}
    \caption{Parameter values used in the simulations carried out for the four-well Hamiltonian in Eq. \eqref{fourwells_ham}.} 
    \label{table_fourwell}
\end{table}
For each of the cases mentioned above, we have taken a random sample of $10^4$ initial conditions on the Poincar\'e surface of section $y = 0$, $p_y \geq 0$. This process has provided us with a large dataset to test how well our trained SVM model performs when detecting chaos and regularity in Hamiltonian systems defined by different types of PES.

Finally, the Hamiltonian model that we have used to test the performance of our SVM model is the classical H\'{e}non-Heiles system \cite{henon1964applicability}, which was introduced by Michel H\'enon and Carl Heiles in 1964 as a model to understand the motion of stars around a galactic center. This Hamiltonian is a prototypical example used to understand the fundamental concepts of chaos, the structure of phase space, and the impact of energy on dynamical stability. From this system, one can easily explore the transition from regular to chaotic motion, and how the onset of chaos depends on the energy of the system. It has been shown that for low energies, this system displays regular and quasi-periodic, while above a certain threshold, the motion becomes chaotic. The H\'enon-Heiles Hamiltonian is defined by:
\begin{equation}
    \mathcal{H}(x,y,p_x,p_y) = \dfrac{p_x^2}{2} + \dfrac{p_y^2}{2} + V(x,y) \;.
    \label{Ham_HH}
\end{equation}
where the potential energy is:
\begin{equation}
    V(x,y) = \dfrac{1}{2}(x^2+y^2) + x^2y - \dfrac{1}{3}y^3
    \label{pes_HenonHeiles}
\end{equation}
This PES has three symmetrically located index-1 saddles, and a local minimum at the origin, as shown in Fig. \ref{PESs} C). The equations of motion for this system are:
\begin{equation}
\begin{cases}
\dot{x} = \dfrac{\partial H}{\partial p_x} = p_x \\[.4cm]
\dot{y} = \dfrac{\partial H}{\partial p_y} = p_y \\[.4cm]
\dot{p}_x = -\dfrac{\partial \mathcal{H}}{\partial x} = - x - 2xy \\[.4cm]
\dot{p}_y = -\dfrac{\partial \mathcal{H}}{\partial y} = - y - x^2 + y^2 
\end{cases}
\;.
\end{equation}
In this work we have considered the energies $\mathcal{H} = 1/20,1/15,1/12,1/10,1/8$ and, for each of these cases, we have taken a random sample of $10^4$ initial conditions on the Poincar\'e surface of section $x = 0$, $p_x \geq 0$.

\end{appendices}

\bibliography{referencias}

\end{document}